\documentclass[a4paper,times,3p]{elsarticle}
\usepackage[utf8]{inputenc}

\usepackage[T1]{fontenc}

\usepackage{microtype}

\usepackage{booktabs}

\usepackage{gnuplot-lua-tikz}

\usepackage{amsmath}
\usepackage{amssymb}
\usepackage{amsthm}
\usepackage{bm}
\usepackage{graphicx}
\usepackage{subfigure}
\usepackage{epstopdf}
\usepackage{psfrag}
\usepackage{color}

\usetikzlibrary{arrows.meta}

\usepackage{mathtools}

\usepackage[font=footnotesize,labelfont=bf]{caption}
\usepackage{multirow}

\usepackage{pxfonts}

\usepackage{import}

\graphicspath{./graphics/}

\usepackage[textsize=scriptsize]{todonotes}
\usepackage{setspace}

\usepackage{hyperref}

\journal{Computer Methods in Applied Mechanics and Engineering}

\def\statement{\begin{minipage}[t]{.75\textwidth}
       NOTICE: This is the author's version of a work that was accepted for
publication in Computer Methods in Applied Mechanics and Engineering.  Changes
resulting from the publishing process, such as editing, corrections, structural
formatting, and other quality control mechanisms may not be reflected in this
document. Changes may have been made to this work since it was submitted for
publication.
\\

\copyright \, 2017. This manuscript version is made available under the
CC-BY-NC-ND 4.0 license \\
\url{http://creativecommons.org/licenses/by-nc-nd/4.0/}.
       \end{minipage}}

\makeatletter
\def\ps@pprintTitle{%
     \let\@oddhead\@empty
     \let\@evenhead\@empty
     \def\@oddfoot{\footnotesize\itshape
       \statement\hfill\today}%
     \let\@evenfoot\@oddfoot}
\makeatother

\begin{document}
\begin{frontmatter}

\title{Boundary-Conforming Free-Surface Flow Computations: Interface Tracking
for Linear, Higher-Order and Isogeometric Finite Elements}
\author{Florian Zwicke\corref{cor1}}
\ead{zwicke@cats.rwth-aachen.de}
\author{Sebastian Eusterholz\corref{}}
\ead{eusterholz@cats.rwth-aachen.de}
\author{Stefanie Elgeti\corref{}}
\ead{elgeti@cats.rwth-aachen.de}

\cortext[cor1]{Corresponding author}

\address{Chair for Computational Analysis of Technical Systems (CATS), RWTH
Aachen University, 52056 Aachen, Germany\\
Center for Computational Engineering Science (CCES), RWTH Aachen University,
52056 Aachen, Germany}

\begin{abstract}
The simulation of certain flow problems requires a means for modeling a free
fluid surface; examples being viscoelastic die swell or fluid sloshing in tanks.
In a finite-element context, this type of problem can, among many other options,
be dealt with using an interface-tracking approach with the
Deforming-Spatial-Domain/Stabilized-Space-Time (DSD/SST) formulation. A
difficult issue that is connected with this type of approach is the
determination of a suitable coupling mechanism between the fluid velocity at the
boundary and the displacement of the boundary mesh nodes. In order to avoid
large mesh distortions, one goal is to keep the nodal movements as small as
possible; but of course still compliant with the no-penetration boundary
condition. Standard displacement techniques are full velocity, velocity in a
specific coordinate direction, and velocity in normal direction.  In this work,
we investigate how the interface-tracking approach can be combined with
isogeometric analysis for the spatial discretization. If NURBS basis functions
of sufficient order are used for both the geometry and the solution, both a
continuous normal vector as well as the velocity are available on the entire
boundary. This circumstance allows the weak imposition of the no-penetration
boundary condition. We compare this option with an alternative that relies on
strong imposition at discrete points. Furthermore, we examine several coupling
methods between the fluid equations, boundary conditions, and equations for the
adjustment of interior control point positions.
\end{abstract}

\end{frontmatter}

\section{Introduction} \label{sec-intro}

This paper can be placed in the field of free boundary problems. More
specifically, it considers fluid flow problems where the computational domain is
part of the solution --- as in computational domains that contain a free
surface. Examples considered here are sloshing tanks --- under a seismic load,
the liquid stored in a tank begins to slosh --- and die swell --- when a
material melt exits the shape-giving die through which it has been pushed, it
will expand.

Within a finite-element context, different methods are available to follow the
position of the free surface throughout the computation \cite{Elgeti2015,
Caboussat2005}. We generally distinguish between interface-capturing and
interface-tracking methods. Examples for interface capturing methods are
particle methods \cite{Easton72, Girault76}, level-set \cite{Sethian99b,
Osher2001}, volume of fluid \cite{Nichols71, Hirt75, Noh76, Hirt81a}, or phase
field \cite{Gibbs1873, Onsager1931, Prigogine1966, Gyrmati1970, deGroot1984,
Emmerich2002}. These approaches have in common that they use an additional
function/marker to indicate the position of the free surface on an Eulerian
background mesh. Examples for interface tracking are full Lagrangian or
Arbitrary Lagrangian Eulerian (ALE) methods \cite{Hirt74a}. Moreover, we would
like to mention the Deforming-Spatial-Domain/Stabilized-Space-Time (DSD/SST)
formulation \cite{Tezduyar92a}, the interface tracking approach which has been
used in this work.

For the considered test case of sloshing tanks, most methods mentioned above
have been used in literature; as reviewed in \cite{Elgeti2015}.  In more recent
publications, Grotle et al. employ the level-set method to compute turbulent
sloshing \cite{Grotle2016}. A particle method is the basis of the analysis of a
sloshing tank with two compartments performed in \cite{Manderbacka2016}.
The die swell test case is an example with moderate deformations and no
topological changes. It is therefore an ideal candidate for an interface
tracking approach \cite{Elgeti2010, Pauli2012, Stavrev2015, Knechtges2014,
Knechtges2015}. Nevertheless, also other methods are applicable. These include
particle methods \cite{Tome2008, Tome2012, Mompean2011} and volume of fluid
\cite{Rocco2010}.

One major advancement in, but not limited to, the field of finite element
methods is Isogeometric Analysis (IGA). Originally presented in
\cite{Hughes2005}, numerous works have been published on the subject in the last
decade. The key idea of IGA is to replace the Lagrangian shape functions, which
are generally used by finite element frameworks, by the geometry definition, as
it would be used in Computer-Aided-Design (CAD) software. As detailed in
\cite{Cottrell2009}, key advantages of the IGA approach are: (1) exact (as in
CAD-conforming) geometry representation, (2) user-controlled smoothness and
continuity of the basis, (3) in general higher accuracy per degree of freedom,
(4) in general no node-to-node oscillations, etc.
In the broad area of free boundary problems, Isogeometric Analysis has been
applied in the context of the phase field method for the Cahn-Hilliard equations
\cite{Gomez2008}, brittle fracture \cite{Borden2012}, and topology optimization
\cite{Dede2012}. In connection with level-set, IGA has been utilized to compute
the dam break problem \cite{Akkerman2011,Amini2014} or as a boundary indicator
within the finite cell method \cite{Rank2012}. Another free-surface-related
application is the computation of wave resistance on a ship hull using the
isogeometric boundary element method \cite{Ginnis2014}. To the knowledge of the
authors, IGA has not been utilized in conjunction with interface tracking of
free surfaces.

The aim of this paper is therefore to define appropriate mesh displacement
techniques for interface tracking of a free fluid boundary in an Isogeometric
Analysis context. As compared to the classic, well-established methods in
standard finite elements --- displacement with the full fluid velocity, the
normal fluid velocity, and the velocity in one specific coordinate direction ---
there will be one key difference: the unknown velocities and displacements no
longer have a direct representation on the free surface. The reason is that, in
IGA, the geometry is represented by non-uniform rational B-splines (NURBS) or
other spline types. Their shape is controlled via control points. Usually, these
control points will not coincide with the actual geometry.

The paper is structured as follows. Section~\ref{sec-equa} will introduce the
governing equations for fluid flow, the no-penetration boundary condition for
the free surface and the governing equations for the mesh adaptation. Classic
numerical methods for their solution are reviewed in Section~\ref{sec-meth}.
Section~\ref{sec-impl} will then discuss how these equations can be interpreted
in an IGA context. Two numerical examples --- sloshing tank and die swell, both
in 2D --- will be presented in Section~\ref{sec-test}.

\section{Governing equations for the free-boundary value problem of free-surface
flow} \label{sec-equa}


For free-surface flow, three types of equations are relevant: (1) the
Navier-Stokes equations in combination with an appropriate constitutive equation
govern the fluid flow, (2) the displacement of the free surface, which is
governed by the no-penetration boundary condition, and (3) the equations
governing the possible displacement of interior parts of the domain in order to
maintain validity of the mesh. This section will introduce those three equations
in the form in which they have been used for the test cases shown in
Section~\ref{sec-test}.

\subsection{Governing equations for fluid flow: The Navier-Stokes equations}

In the generic incompressible and isothermal fluid flow problem the unknowns are
the velocity, \( {\bf u}({\bf x},t) \), and the pressure, \( p({\bf x},t) \).
The computational domain at each instant in time, denoted by \(\Omega_t \), is a
subset of \(\mathbb{R}^{nsd}\), where \(nsd\) is the number of space dimensions.
Then at each point in time \(t \in [0,T] \), the flow problem is governed by the
Navier-Stokes equations, which in our notation read:

\begin{alignat}{3}
  \label{NavierStokes1}
  \rho \left( \frac{\partial {\bf u}}{\partial t}
    + {\bf u \cdot \nabla u} - {\bf
    f} \right) -  \nabla \cdot \bm{\sigma} &= {\bf 0}
    \quad &&\mbox{on} \quad \Omega_t ,\ t \in [0,T] \,, \\
  \label{NavierStokes2}
  {\bf \nabla \cdot u } &= 0
    \quad &&\mbox{on} \quad \Omega_t ,\ t \in [0,T] \,,
\end{alignat}
with \( \rho\) as the density of the fluid. We consider only Newtonian fluids,
meaning that the stress tensor \( \bm{\sigma}\) is defined as

\begin{align}
\bm{\sigma} ( {\bf u}, p) &= -p {\bf I} + 2 \mu \bm{\varepsilon}( {\bf u}) \,,
\label{eqn:Newtonian}
\end{align}
with
\begin{align} \label{eqn:strain}
\bm{\varepsilon}({\bf u}) &= \frac{1}{2} (\nabla {\bf u} + (\nabla {\bf u})^T) \,,
\end{align}
where \( \mu  \) denotes the dynamic viscosity.  $ {\bf f} $ includes all
external body forces with respect to the unit mass of fluid. Note that the
spatial domain is time-dependent, which is indicated by subscript \(t\).

For creeping flows (i.e., Reynolds number $\ll1$), the advective term in the
Navier-Stokes equations is often neglected, giving rise to the Stokes equations:

\begin{alignat}{3}
  \label{StatStokes_1}
  \rho \left( \frac{\partial {\bf u}}{\partial t}
    - {\bf f} \right) -  \nabla \cdot \bm{\sigma} &= {\bf 0}
    \quad &&\mbox{on} \quad \Omega_t ,\ t \in [0,T] \,, \\
  \label{StatStokes_2}
  {\bf \nabla \cdot u } &= 0
    \quad &&\mbox{on} \quad \Omega_t ,\ t \in [0,T] \,.
\end{alignat}
The constitutive equations for Newtonian fluids remain the same as above.

In the transient case, a divergence-free velocity field for the whole
computational domain is needed as an initial condition:

\begin{align}
  {\bf u} ( {\bf x}, 0 ) = \hat{{\bf u}}^0( {\bf x})
  \quad \mbox{on} \quad \Omega_t ,\ t=0 \,.
\end{align}

In order to obtain a well-posed system, boundary  conditions have to be imposed
on the external boundary of $\Omega_t$, denoted as $\Gamma_t$. Here, we
distinguish between Dirichlet and Neumann boundary conditions given by:

\begin{alignat}{3}
  {\bf u} &= \hat{{\bf u}}
    \quad &&\mbox{on} \quad \Gamma_u,\ t \in [0,T] \,, \\
  {\bf n} \cdot \bm{\sigma} &= \hat{{\bf h}}
    \quad &&\mbox{on} \quad \Gamma_h,\ t \in [0,T] \label{eqn:Neumanncondition} \,,
\end{alignat}

\noindent where $\hat{{\bf u}}$ and $\hat{{\bf h}}$ are prescribed velocity and
stress values. $ \Gamma_u $ and $ \Gamma_h $ denote the Dirichlet and Neumann
part of the boundary, forming a complementary subset of $\Gamma_t$, i.e., $
\Gamma_u \cup \Gamma_h = \Gamma_t$ and $ \Gamma_u \cap \Gamma_h = \emptyset$.
{\bf n} refers to the outer normal vector on $ \Gamma_h $.


\subsection{Displacement of the free-surface: The no-penetration boundary condition}
\label{sec:no-pen}
Depending on the specific application the full boundary  $ \Gamma $, or a
portion thereof, can be defined as a free surface $ \Gamma_{free} $. In an
interface-tracking context, we want to ensure that the computational domain
adapts to the fluid flow. The deformation of the free surface is governed by a
no-penetration boundary condition, meaning that no fluid is allowed to cross the
boundary. This condition will be fulfilled as long as the mass flux $\dot{m}$
through  $\Gamma_{free} $ --- and any given subset of  $\Gamma_{free} $ --- is
$0$. The mass flux $\dot{m}$  can be computed as:

\begin{align}
\dot{m} = \int_{\Gamma^*} \rho( {\bf u} - {\bf v} ) {\bf n}\;\text{d}{\bf x}
\overset{!}{=} 0 \,, \quad \forall \, \Gamma^* \subseteq \Gamma_{free}
\label{eqn:no-pen}
\end{align}

with $\rho$ as fluid density, ${\bf u}$ as fluid velocity and ${\bf v}$ as
mesh velocity.

Consequently, all choices for ${\bf v}$ that comply with the kinematic boundary
condition

\begin{align}
{\bf v}({\bf x}) \cdot {\bf n}({\bf x}) = {\bf u}({\bf x}) \cdot {\bf n}({\bf x})
\label{eqn:kbc}
\end{align}

\noindent are valid.

The straightforward choice is

\begin{align}
  {\bf v}({\bf x}) = {\bf u} ({\bf x}) \,.
\label{eqn:fullLagrange}
\end{align}

This may however not always be the best choice. Equation~\eqref{eqn:kbc}
indicates that the tangential component of  {\bf v} has no influence on the
boundary shape --- which is the actual result of the boundary displacement ---,
but may significantly influence the quality of the utilized computational mesh.
Particularly applications with large tangential velocities --- such as die swell
---  profit tremendously if the tangential velocity component is modified or
even suppressed when computing the mesh velocity ${\bf v}$.


Behr \cite{Behr92} indicates several possibilities for boundary displacement
resulting from Equation~\eqref{eqn:kbc}: displacement with the normal velocity
component with $ {\bf v} = ({\bf u} \cdot {\bf n}){\bf n} $ and displacement
only in a specific direction ${\bf d}$ (e.g., $y$-direction), i.e., $ {\bf v} =
\frac{({\bf u} \cdot {\bf n}) {\bf d}}{ {\bf n} \cdot {\bf d}} $. The concept is
detailed in \cite{Elgeti2015}.

\subsection{Retaining mesh quality: Mesh update for inner nodes}

As the boundary  $ \Gamma_{free} $ is modified, usually the discretization of
the interior of the domain needs to be adapted as well. For this purpose, we
employ the Elastic Mesh Update Method (EMUM) \cite{Johnson94a}. In this method,
the computational mesh is treated as an elastic body reacting to the boundary
deformation applied to it. The linear elasticity equation is solved for the mesh
displacement ${\bf z}$:

\begin{align}
  \label{eq:sopt-emum}
  \nabla \cdot \mbox{${\bm{\sigma}}$}_{\mathrm{mesh}} &= {\bf 0} \,, \\
  \mbox{${\bm{\sigma}}$}_{\mathrm{mesh}}({\bf z}) &= \lambda_{mesh} \left(
    {\mathrm{tr}} \, {\bm{\varepsilon}}_{\mathrm{mesh}} ({\bf z})\right){\bf
    I} + 2 \mu_{mesh} {\bm{\varepsilon}}_{\mathrm{mesh}}({\bf z}) \,, \\
  {\bm{\varepsilon}}_{\mathrm{mesh}} ({\bf z}) &= \frac{1}{2} \left(\nabla
    {\bf z} + (\nabla {\bf z})^T\right) \,.
\end{align}
 \( \lambda_{mesh} \) and \( \mu_{mesh} \) --- in structural mechanics the
Lam\'e-parameters --- have no physical meaning within the mesh deformation. They
can be chosen freely for each element in order to control its respective
stiffness.

The mesh update equations are solved at discrete points in time $t$ on the
domain $\Omega_t$. We partition the boundary $\Gamma$ into the disjoint parts
$\Gamma_{free}$, $\Gamma_{fixed}$ and $\Gamma_{slip}$ with
$\Gamma = \Gamma_{free} \cup \Gamma_{fixed} \cup \Gamma_{slip}$.
$\Gamma_{fixed}$ is the part of the boundary that does not change at all
and $\Gamma_{slip}$ is the part that allows tangential node movement.
We obtain the boundary conditions

\begin{alignat}{2}
  {\bf z} &= {\bf 0} &&\quad \mbox{on}
    \quad \Gamma_{fixed} \,, \\
  {\bf z} &= {\bf z}_D &&\quad \mbox{on}
    \quad \Gamma_{free} \,, \\
  {\bf z} \cdot {\bf n} &= 0 &&\quad \mbox{on}
    \quad \Gamma_{slip} \,, \\
  \left(\mbox{${\bm{\sigma}}$}_{\mathrm{mesh}}({\bf z}) {\bf n}\right)
    \cdot {\bf t}_l &= 0
    &&\quad \mbox{on} \quad \Gamma_{slip} \,,
    \ 1 \le l \le d-1 \,,
\end{alignat}

where ${\bf z}_D$ is the displacement on the free-surface boundary, ${\bf n}$ is
the surface normal vector and ${\bf t}_l$ are a set of $d-1$ linearly independent
surface tangent vectors.

One should note that the slip boundary condition as it is described here can
only preserve the shape of the boundary if it is completely straight.

\section{Background on the employed simulation methods} \label{sec-meth}

\subsection{The Deforming-Spatial-Domain/Stabilized-Space-Time (DSD/SST)
formulation}

\begin{figure}[h]
\center
\psfrag{a}{free surface at}
\psfrag{b}{time step $t^{n+1}$}
\psfrag{c}{space-time}
\psfrag{d}{element}
\psfrag{e}{free surface at}
\psfrag{f}{time step $t^{n}$}
\includegraphics[width=10cm]{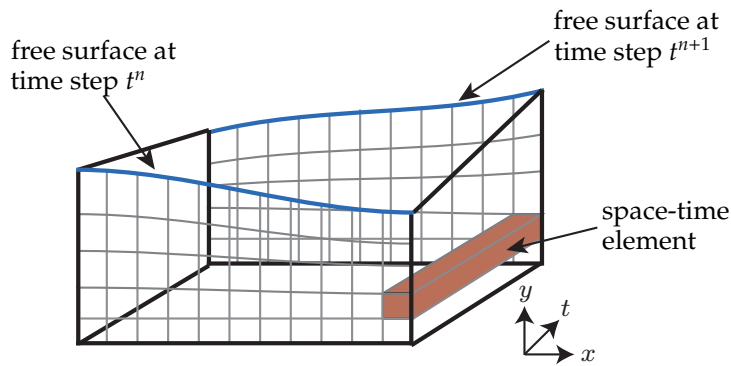}
\caption{Illustration of a space-time slab tracking the deformation of the
domain. In the front, we see the old time step, in back the new time step with a
deformed domain.} \label{fig:space-time}
\end{figure}

The method we use for simulating our test cases is a variant of the
deforming-spatial-domain/stabilized-space-time finite element formulation. As
the name already indicates, it is a finite element method that is based on a
space-time discretization --- i.e., both the spatial domain and the temporal
interval are discretized using finite elements. Similar to time-stepping, the
space-time domain is usually subdivided into so-called space-time slabs (Figure
\ref{fig:space-time}). It is the use of space-time finite elements that is
responsible for the main advantage of the method: the variational form of the
governing equations can be written over a deforming domain.
This allows mesh displacement in a very natural way, since the governing
equations do not require any adjustments.
Specifically, in contrast to standard ALE methods, the mesh velocity does not
need to be included into the formulation as an advection velocity.
Details on the method
can be found in \cite{Tezduyar92a,Tezduyar92b}. Other extensions of DSD/SST are
for example used by Tezduyar~\cite{Tezduyar2004}, Takizawa et al.
\cite{Takizawa2011, Takizawa2014} and  Zilian et al. \cite{Zilian2008}.
Since the Navier-Stokes equations require stabilization when they are solved
with the finite element method, we add Galerkin/Least-Squares (GLS) terms to the
discretized weak form.

\subsection{Displacement of the free surface for standard finite elements}

\begin{figure}[h]
\center
\psfrag{a}{possible displacement velocity of an FE node}
\psfrag{b}{line of possible displacement velocities}
\psfrag{c}{for the FE node}
\includegraphics[trim = -15mm 0mm 0mm 0mm, clip,
width=10.2cm]{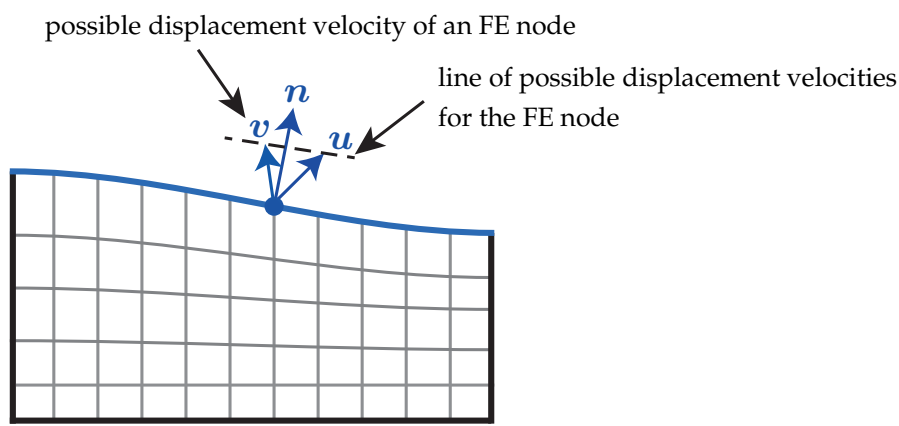}
\caption{Illustration of possible displacement velocities for a given FE node on
the free surface obeying the no-penetration boundary condition.}
\label{fig:FEdisp}
\end{figure}

As discussed in Section \ref{sec:no-pen}, each point on the surface has to be
displaced according to the rule that the normal component of the displacement
velocity is equal to the normal component of the fluid velocity (cf.
Equation~\eqref{eqn:kbc}). For standard finite element discretizations, the
finite element nodes lie directly on the surface; even more, they constitute the
approximation of the surface. One approach to displacing the free surface is
therefore to apply the chosen displacement rule (e.g., normal direction,
vertical direction) directly to the finite element nodes. The procedure is
illustrated in Figure~\ref{fig:FEdisp}. The challenge is the computation of the
normal vector. Since the standard finite element shape functions are
$C^0$-continuous across element boundaries, that vector is not directly
definable via the finite element approximation, but usually computed as a
weighted average of the normals of surrounding element faces \cite{Engelman82a}.
Other options are the computation from additional boundary descriptions that
complement the finite element approximation \cite{Elgeti2010,Stavrev2015}.

It is to be noted that this procedure will result in only an approximation of
the no-penetration boundary condition (Equation~\eqref{eqn:no-pen}) as the
displacement rule is not fulfilled exactly in between nodes (the normal vector
on a face is no longer constructable via the normal vectors at the corresponding
nodes); however, an approximation is inevitable, as we only have an
approximation of our geometry available in any case.

In summary, note that this method is tailored to linear temporal and spatial
discretizations. It will perform poorly for higher-order shape functions. A
requirement is that the nodes lie on the boundary. The accuracy of the method
depends strongly on the type of approximation for the normal vector at the
finite element nodes.

\subsection{Isogeometric Analysis}

Modern engineering design evolves around \textit{Computer-Aided Design} (CAD)
systems. As compared to drawing boards, modern CAD offers a significantly
enhanced experience to the user. In engineering design, computer based geometry
description is almost exclusively based on splines. More specifically, most CAD
systems use Non-Uniform Rational B-Splines (NURBS) as a means of creating,
communicating, and storing geometry information. From a strict engineering point
of view, e.g., in a CAD/CAM context, this representation has proven to be very
efficient.

\subsubsection{Non-uniform rational B-splines (NURBS)}

Detailed descriptions of NURBS and their properties can for example be found in
\cite{Piegl97,Rogers2001}. These books are also the basis for the following
short description. NURBS belong to the category of parametric geometry
representations. Their concept is derived for 1D geometrical entities, but can
be extended to any spatial dimension. One-dimensional and two-dimensional
splines are most common, however. In many cases, we encounter scenarios, where
lower-dimensional spline entities are immersed in a higher-dimensional space. A
NURBS curve ${\bf C}$ (i.e., the 1D entity) is defined as

\begin{align}
{\bf C}(\theta) = \sum_{i=1}^n R_{i,p}(\theta) {\bf P}_i,
\label{eqn:NURBScurve}
\end{align}

\noindent where $R_{i,p}$ denote the NURBS basis functions ($p$ being the
polynomial degree) and $ {\bf P}_i$ are the coordinates of the control point $i$
\--- which are always given in the global coordinate system ${\bf x}$. Note
again that the dimension of ${\bf x}$ can be $1$ but also larger. $n$ indicates
the total number of control points. In the following, we will discuss the
constituents of a NURBS curve in more detail.

The NURBS basis functions $R_{i,p}$ are rational polynomial functions with local
support. The beginning and end of the interval, where a certain basis function
is nonzero, is marked by a so-called knot. Knots are simply values in the
parameter space. The range of $R_{i,p}$ is $[0,1]$, however, the value of $1$ is
only rarely reached. One implication of this property is that the curve usually
does not coincide with the control points.

The control points $ {\bf P}_i$ define the position and shape of the NURBS
entity. They are arranged in a structured grid. Depending on the support of the
corresponding basis function, which spans at most $p+1$ elements, each control
point influences the shape of the NURBS entity only in a limited area.

One important aspect of NURBS are their continuity properties. Between the basis
and the distribution of control points, the user can control the continuity
almost arbitrarily; constructing kinks or other discontinuities when needed, but
profiting from high continuity otherwise.
Specifically, a NURBS of polynomial degree $p$ will be $p-1$ times continuously
differentiable when none of the knot values are repeated.

\subsubsection{The isoparametric concept}

Cottrell, Bazilevs and Hughes used NURBS to devise the concept of isogeometric
analysis (IGA) \cite{Hughes2005,Cottrell2009,Bazilevs2010}. In the spirit of the
isoparametric concept, IGA utilizes the NURBS basis functions in order to
represent both the geometry and the unknown solution. For the geometry, the
already presented definition \eqref{eqn:NURBScurve} is employed. The unknown
function ${\bf u}$ is then approximated in exactly the same fashion as:

\begin{align}
{\bf u}^h = \sum_{i=1}^n R_{i,p}(\theta) {\bf u}_i \,.
\end{align}

\begin{figure}
\psfrag{a}{$u_1$}
\psfrag{b}{$u_2$}
\psfrag{n}{$u_3$}
\psfrag{c}{$u_4$}
\psfrag{d}{$u_5$}
\psfrag{e}{$u^h$}
\psfrag{f}{$x$}
\psfrag{g}{$u^h(x)$}
\psfrag{h}{${\bf P}_1$}
\psfrag{i}{${\bf P}_2$}
\psfrag{m}{${\bf P}_3$}
\psfrag{j}{${\bf P}_4$}
\psfrag{k}{${\bf P}_5$}
\psfrag{l}{$\theta$}
\psfrag{x}{$x$}
\centering
\includegraphics[scale=0.3]{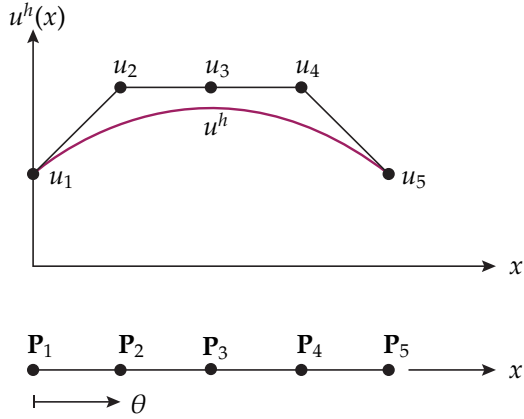}
\caption{Approximation of an unknown function $u$ on an isogeometric grid: The
lower part of the picture shows a rod, which is represented using a quadratic
NURBS with $5$ control points. The upper part of the picture shows the
representation of the unknown function $u$ related to this grid. Instead of the
five control points for the geometry, we now have five control variables that
define the solution $u^h$.}
\label{fig:iga-function-representation}
\end{figure}

The values ${\bf u}_i$ are the control variables. These are the unknown values
solved for in the finite element code. Note that just like the control points
for the geometry, the values of the control variables are in general not
coinciding with the solution; see Figure~\ref{fig:iga-function-representation}.

This idea not only essentially provides the capability to perform finite element
analysis directly on CAD models, but it also profits from superior approximation
properties of NURBS as compared to the polynomial shape functions used in the
standard finite element method.

When used in free-surface flow, IGA is traditionally combined with an interface
capturing approach. Recall that interface tracking, in contrast, means that the
mesh will deform with the free surface. When IGA is applied to interface
tracking, the challenge will be to find a proper displacement rule for the
control points, which describe the free surface. Unlike in standard finite
elements, the control points, which define the shape of the NURBS, are not
located on the surface. The kinematic boundary condition
(Equation~\eqref{eqn:kbc}), however, needs to be fulfilled exactly there. The
question will be, how to displace the control points in such a way, that the
resulting curve displacement complies with the kinematic boundary condition?

\begin{figure}[h]
\center
\psfrag{a}{normal vector available}
\psfrag{b}{at the free surface}
\psfrag{c}{fluid velocity available}
\psfrag{d}{at control point}
\includegraphics[trim = -15mm 0mm 0mm 0mm, clip,
width=10.0cm]{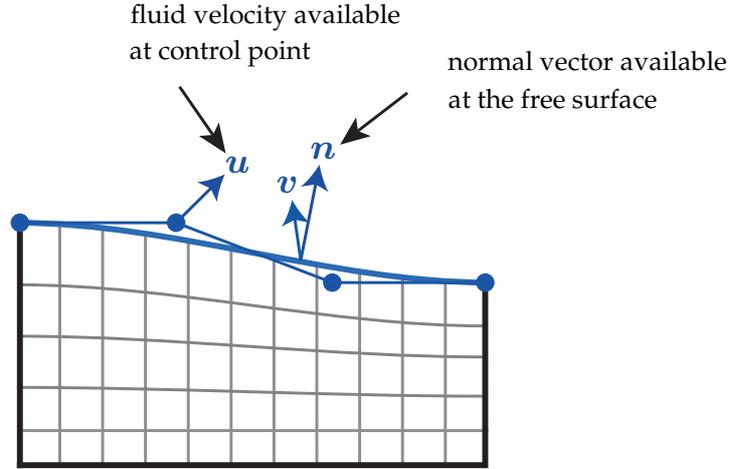}
\caption{Illustration of possible displacement velocities for a given point on
the free surface obeying the no-penetration boundary condition.}
\label{fig:IGAdisp}
\end{figure}

\section{Implementation: How should the free surface be displaced?} \label{sec-impl}

As described in Section~\ref{sec:no-pen}, the displacement of the free surface
needs to fulfill the no-penetration boundary condition
(Equation~\eqref{eqn:no-pen}), which can be rewritten into the kinematic
boundary condition (Equation~\eqref{eqn:kbc}). Even though setting the mesh
displacement velocity equal to the fluid velocity will always fulfill this
condition, it is vital to most application cases that other displacement schemes
--- e.g., in normal or in vertical direction --- can be employed. Otherwise,
tangential velocity components will have a negative influence on the mesh
quality.

When this displacement concept is transferred to IGA, we face the problem that
the unknowns (fluid and mesh velocity) are connected to the control points,
whereas the kinematic boundary condition needs to hold for every point {\bf on}
the spline. The same holds for the normal vector, which can only be computed for
the actual spline. At the control points, we have no notion of a normal vector.

Due to similar issues, weak imposition of Dirichlet boundary conditions has
become a serious option for IGA \cite{Bazilevs2007}. Therefore, in addition to
variants of the kinematic boundary condition (Equation~\eqref{eqn:kbc}) ---
which in a sense can be considered a strong imposition of a boundary condition
--- methods based on the integral formulation of the no-penetration boundary
condition (Equation~\eqref{eqn:no-pen}) --- showing similarities to a  weak
imposition --- might be beneficial.

In this paper, we propose two displacement methods for the free surface; one
based on strong imposition and one on weak imposition. The first option entails
moving the control points with normal vectors that are computed at the Greville
abscissae. The second option considers the weak fulfillment of the
no-penetration boundary condition --- meaning that an additional equation system
will need to be solved --- in three different variants: (1) displacement in the
surface normal direction, (2) displacement with the full velocity, and (3)
displacement in a certain fixed direction.

\subsection{Displacement based on the normal vector at the Greville abscissae}

As detailed before, there is no notion of a normal vector at a control point.
One point on the spline, whose normal vector might come close to a vector that
might be considered a normal vector at a specific control point, is the normal
vector at the corresponding Greville abscissa. The Greville abscissa is the
point on the surface the control point converges to in case of refinement. The
Greville abscissa can be computed as the average of the knot values relevant for
a control point, excluding the first and last value \cite{Farin2014}. Typically,
this is also close to the point where the associated basis function is maximal
(this is exactly fulfilled for uniform knot vectors). The normal vector at any
parametric coordinate $\theta$ of a NURBS curve can then be computed as
\cite{Gray1998}:
 \begin{align}
	\bm{n}(\theta) =\left(\begin{array}{c}-t_{y} \\t_{x}\end{array}\right),
	\qquad\text{with tangent vector }\bm{t}(\theta)=\left(\begin{array}{c}t_{x} \\t_{y}\end{array}\right) 
	= \frac{\bm{C}'(\theta)}{\left|\bm{C}'(\theta)\right|}.
\end{align} 

The formula requires $\bm{C}'(\theta)$, the first derivative of the curve with
respect to the local parameter $\theta$, for the definition of which we refer to
\cite{Piegl97}.

Based on the fluid velocity at any given control point ${\bf P}_i$, we define
its displacement (or control point coordinate increment) $\Delta {\bf P}_i$ as:

\begin{align}
\Delta {\bf P}_i \coloneqq \left( \bar{{\bf u}}_i \cdot \bar{{\bf n}}_i \right)
\bar{{\bf n}}_i \cdot \Delta t \,,
\end{align}

with $\bar{{\bf u}}_i$ and $\bar{{\bf n}}_i$ being the fluid velocity at control
point $i$ and the normal vector evaluated at the $i^{th}$ Greville abscissa,
both averaged over the time slab.

\subsection{Displacement based on weak formulation of the kinematic boundary
condition}

We can also treat the no-penetration condition as a partial differential
equation (PDE) on the boundary. We can solve this with the finite element
method. Since the no-penetration condition provides just one scalar equation, it
is not sufficient to solve for multiple boundary coordinate values in two- or
three-dimensional space.

Multiple options are available to formulate a well-determined system of
equations. To cover the different options, we will first treat the
discretization of the equations in an abstract way by considering an operator
$\mathcal{F}({\bf u}, {\bf s})$ and the PDE

\begin{align}
  \mathcal{F}({\bf u}, {\bf s}) &= 0
  \quad \text{on} \quad \Gamma_{free,t} \, , \, t \in [0, T] \, .
  \label{eqn:fpde}
\end{align}

In this equation, ${\bf u}$ is the fluid velocity and ${\bf s}$ is the boundary
displacement relative to a reference configuration $\Omega_0$.
We denote the free-surface portion of the boundary of $\Omega_t$ as
$\Gamma_{free,t}$.

We introduce a spatial discretization space
$\mathcal{I}_0 \subset C^0\left(\Gamma_{free,0}\right)$.
This is a finite-dimensional space that may be based on Lagrange or spline
interpolation, depending on whether the standard finite element method or
isogeometric analysis is used.
We look at a single time slab $J = [t_k, t_{k+1}] \subseteq [0, T]$.
We define a function

\begin{align}
  &\tilde{\varphi} \in \mathcal{I}_0^d \otimes \mathbb{P}_1(J) \, ,
\end{align}

with spatial dimensionality $d$, such that

\begin{align}
  &\tilde{\varphi} : \, \Gamma_{free,0} \times J \to \tilde{P}
  \quad \text{with} \quad \tilde{P} \subset \mathbb{R}^d \, .
\end{align}

We now choose $\tilde{\varphi}$ such that $\varphi$, defined by

\begin{align}
  \varphi({\bf x}, t) = \left( \begin{array}{c}
    \tilde{\varphi}({\bf x}, t) \\ \text{Id}_J(t)
  \end{array} \right) \, , \\
  \text{with} \quad
  \varphi : \, \Gamma_{free,0} \times J \to P_{free} \, ,
\end{align}

is a homeomorphism.
Here, $P_{free}$ is the deformed space-time
free-surface boundary, which can be related to $\Gamma_{free,t}$ by

\begin{align}
  P_{free} &= \left\{ ({\bf x}, t) \in \mathbb{R}^d \times J
              \, \middle| \, {\bf x} \in \Gamma_{free,t} \right\} \, .
\end{align}

We use this to define a trial space $\mathcal{S}$ and test space $\mathcal{T}$:

\begin{align}
  \mathcal{S} &= \left\{ s \circ \varphi^{-1} \, \middle| \,
                 s \in \mathcal{I}_0 \otimes \mathbb{P}_1(J) \right\} \, , \\
  \mathcal{T} &= \left\{ w \circ \varphi^{-1} \, \middle| \,
                 w \in \mathcal{I}_0 \otimes \mathbb{P}_0(J) \right\} \, ,
  \label{eqn:testspace}
\end{align}

where $\mathbb{P}_p$ is the space of polynomials of degree $p$.
One should note that different polynomial degrees in time are used for these
spaces. Specifically, the test space $\mathcal{T}$ only includes functions that
are constant in time.

We want the geometry description to be continuous over time slab
boundaries. Therefore, we have to impose a temporal Dirichlet boundary condition
on the lower level of the time slab for the discrete displacement field ${\bf s}^h$.
This requires an adjusted trial space $\mathcal{S}_s$:

\begin{align}
  \mathcal{S}_s &= \left\{ {\bf s} \in \mathcal{S}^d \, \middle| \,
                 \left. {\bf s} \right|_{t_k} = {\bf s}_k \right\} \, ,
\end{align}

where ${\bf s}_k$ with
${\bf s}_k \circ \left. \varphi \right|_{t_k} \in \mathcal{I}^d_0$ is a
given initial condition for the current time slab.
The relationship between the discrete displacement field ${\bf s}^h$ and the
geometry is given by

\begin{align}
  \tilde{\varphi}^{-1}({\bf x}, t) &= {\bf x} - {\bf s}^h({\bf x}, t)
  \quad \text{on} \quad P_{free} \, .
  \label{eqn:srestrictions}
\end{align}

This places a further restriction on ${\bf s}^h$, since only such choices are
valid which allow $\varphi$ to be a homeomorphism.

As mentioned before, the test space $\mathcal{T}$ in (\ref{eqn:testspace}) was
designed to include only functions that are constant in time. In this way, it
differs from the trial space $\mathcal{S}$ and has a lower dimensionality. This
is necessary, since the dimensionality of $\mathcal{S}_s$ is also reduced. With
the selected spaces, we get

\begin{align}
  d \cdot \text{dim} \, \mathcal{T} &= \text{dim} \, \mathcal{S}_s \, ,
\end{align}

which ensures that the number of equations will match the number of unknown
displacement variables.

We can now formulate the discretized weak form of the equation:
\textit{Find ${\bf u}^h \in \mathcal{S}^d$ and ${\bf s}^h \in \mathcal{S}_s$
such that}

\begin{align}
  \int_{J} \int_{\Gamma_{free,t}} w^h
  \mathcal{F}\left({\bf u}^h, {\bf s}^h\right) \, \text{d}{\bf x} \, \text{d}t &=
  {\bf 0}
\end{align}

\textit{for all $w^h \in \mathcal{T}$.}

In order to form a closed system, these equations have to be combined with the
fluid (Navier-Stokes) and mesh update (EMUM) equations.

We have so far kept the descriptions independent of the spatial dimensionality.
However, while both two- and three-dimensional options for the operator
$\mathcal{F}$ are available, we will at this point stay in the two-dimensional
realm.
In this work, we investigate three possible choices for the operator
$\mathcal{F}$.

A simple option is to move the mesh boundary exactly with the fluid
velocity. We are going to call this option \textit{equal movement} in the
remainder.
For this option, we set the mesh velocity as

\begin{alignat}{3}
  && {\bf v} &= {\bf u} \\
  \Leftrightarrow \quad && \frac{\partial {\bf s}}{\partial t} &= {\bf u} \, .
\end{alignat}

We have seen best results with a variant of the equations that is formulated for
a rotated coordinate system that is aligned with the boundary. This means that
the equations are not formulated using horizontal and vertical parts, but using
parts that are normal or tangential to the boundary.
We will make use of a boundary tangential vector ${\bf t}$ and boundary normal
vector ${\bf n}$, which can both be interpreted as functions of the displacement
field ${\bf s}$.
The condition for equal mesh movement, formulated in normal and tangential
directions, now reads

\begin{align}
  \frac{\partial {\bf s}}{\partial t} \cdot {\bf n} &= {\bf u} \cdot {\bf n} \, , \\
  \frac{\partial {\bf s}}{\partial t} \cdot {\bf t} &= {\bf u} \cdot {\bf t} \, ,
\end{align}

which yields the following choice for $\mathcal{F}$:

\begin{align}
  \mathcal{F}_1({\bf u}, {\bf s})
    &= \frac{\partial {\bf s}}{\partial t} \cdot {\bf n}({\bf s})
     - {\bf u} \cdot {\bf n}({\bf s}) \, , \label{eqn:nopen} \\
  \mathcal{F}_2({\bf u}, {\bf s})
    &= \frac{\partial {\bf s}}{\partial t} \cdot {\bf t}({\bf s})
     - {\bf u} \cdot {\bf t}({\bf s}) \, .
\end{align}

For flows with a strong velocity component tangential to the surface --- e.g.,
die swell --- , this option will lead to large mesh distortions:
the mesh can be swept away!

The first equation for the equal movement is exactly the no-penetration
condition that was formulated previously. A logical alternative option for the
mesh movement results if we prevent mesh movement in the tangential direction.
We call this \textit{normal-direction movement}.
The first equation (\ref{eqn:nopen}) stays in place, and the second one
changes to

\begin{align}
  \mathcal{F}_2({\bf u}, {\bf s})
    &= \frac{\partial {\bf s}}{\partial t} \cdot {\bf t}({\bf s})
\end{align}

in this case.

This choice minimizes
${\bf v} = \partial {\bf s}/\partial t$,
which may help in preserving mesh quality.
Furthermore, it places no restrictions on the direction of ${\bf n}$ or ${\bf
u}$.

Finally, we can also restrict the mesh movement to an arbitrary
\textit{direction}. We still use (\ref{eqn:nopen}) to fulfill the no-penetration
condition, but we set the tangential movement according to a direction vector
${\bf d}$:

\begin{align}
  \mathcal{F}_2({\bf u}, {\bf s})
    &= \frac{\partial {\bf s}}{\partial t} \cdot {\bf t}({\bf s})
     - \frac{{\bf u} \cdot {\bf n}({\bf s})}{{\bf d}
       \cdot {\bf n}({\bf s})}
       {\bf d} \cdot {\bf t}({\bf s}) \, .
\end{align}

Here, ${\bf d}$ signifies a selected direction vector. Usually, ${\bf d}$  will
be defined as one of the unit coordinate vectors. However, it is required that
${\bf d}\cdot {\bf n} \neq 0 $, implying that ${\bf d}$ can never be tangential
to the surface. Note that mesh quality will suffer for small values of ${{\bf
d}\cdot {\bf n}}$. In contrast, for ${{\bf d} \approx {\bf n}}$, this definition
of $ {\bf F}$ can preserve mesh quality over long stretches of time. In
particular, this is due to the fact that a return to the initial boundary will
always exactly recreate the initial boundary mesh (i.e., there can be no
tangential sliding of control points over time). This stands in contrast to the
displacement in normal direction.

For the case of pure vertical movement, which we have used in our examples, the
simpler, but equivalent, equation

\begin{align}
  \mathcal{F}_2({\bf u}, {\bf s})
    &= \frac{\partial {\bf s}_1}{\partial t} \,
\end{align}

can be used.

In the following section, the different displacement options will be applied to
the test cases.

\section{Test cases: How do the PDE-based and node-based displacements compare?}
\label{sec-test}

In order to compare and validate the new approaches to surface displacement, we
have selected two common test cases. The first test case --- a sloshing tank
---- features a damped periodic movement with a strong component orthogonal to
the sloshing direction. The second test case --- die swell --- is also
characterized by a movement orthogonal to the flow direction, however, in this
case, the flow velocity features a strong tangential component.

\subsection{Test case sloshing tank}

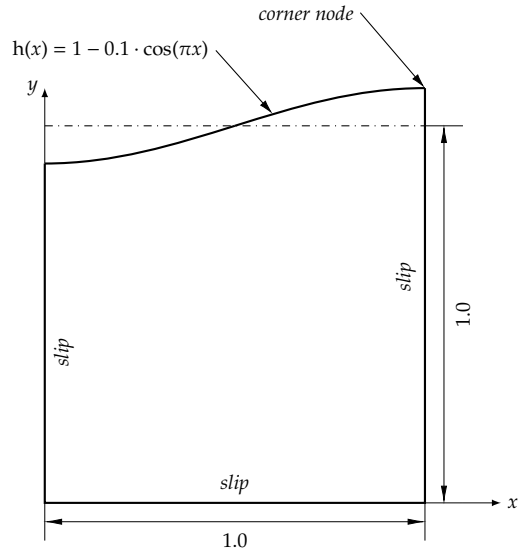
\begin{figure}[h!]
  \centering
  \begin{tikzpicture}[xscale=5, yscale=5]

    \tikzset{every node/.append style={scale=0.75}}

    \draw [thick] (0,0.9) -- (0,0) -- (1,0) -- (1,1.1);
    \draw [thin, dash dot] (0,1) -- (1,1);
    \draw [thick] (0,0.9) cos (0.5,1);
    \draw [thick] (0.5,1) sin (1.0,1.1);

    \draw [very thin] (0.0,0.0) -- (0.0,-0.1);
    \draw [very thin] (1.0,0.0) -- (1.0,-0.1);
    \draw [{Latex[length=2.25mm,width=0.75mm]}-{Latex[length=2.25mm, width=0.75mm]}, very thin] (0.0,-0.05) -- (1.0,-0.05);
    \node [scale=1] at (0.5,-0.1) {$\normalsize 1.0$};
    \draw [very thin] (1.0,0.0) -- (1.1,0.0);
    \draw [very thin] (1.0,1.0) -- (1.1,1.0);
    \draw [{Latex[length=2.25mm,width=0.75mm]}-{Latex[length=2.25mm, width=0.75mm]}, very thin] (1.05,0.0) -- (1.05,1.0);
    \node [rotate=90] at (1.1,0.5) {$1.0$};

    \draw [very thin, >=latex, ->] (0.0,0.9) -- (0.0,1.1) node[left]{$y$};
    \draw [very thin, >=latex, ->] (1.1,0.0) -- (1.2,0.0) node[right]{$x$};
    \node [thick] at (0.5,0.05) {$slip$};
    \node [thick, rotate=90] at (0.05,0.4) {$slip$};
    \node [thick, rotate=90] at (0.95,0.6) {$slip$};
    \draw [-{Latex[length=2mm, width=1mm]}, thin] (0.45,1.2) -- (0.6,1.0309);
    \node [thick, left] at (0.45,1.2) {$\text{h}(x) = 1-0.1 \cdot \text{cos}(\pi x)$};
    \draw [{Latex[length=2mm, width=1mm]}-, thin] (1.0,1.1) -- (0.83,1.3) node[left]{$corner\ node$};

  \end{tikzpicture}
  \caption{Sloshing tank: Setup of geometry, boundary and initial conditions.}
  \label{fig:tank.setup}
\end{figure}

The first test case considers a fluid enclosed in a tank with three straight
walls at a 90-degree angle. As indicated in Figure~\ref{fig:tank.setup}, the
tank has a dimensionless width of $1.0$. The fluid --- viscosity $0.01$ and
density $1000$ --- is initially in a non-equilibrium position. The initial
free-surface position is described via its height $\text{h}(x) = 1-0.1 \cdot
\text{cos}(\pi x)$. The average height is $1.0$. On all walls, a slip boundary
condition is imposed.
This means that the velocity is set to zero in the surface normal direction
using a combination of Dirichlet and Neumann boundary conditions for the
different components.
Upon start of the simulation, the fluid will start to
slosh, heading towards its equilibrium position at constant height. The sloshing
is computed for a dimensionless time interval of $50$. This corresponds to ten
full sloshing cycles given a gravity constant of $-1.0$.

For all standard finite element computations, we use a Q1Q1 (i.e., bilinear basis
functions for both pressure and velocity interpolation spaces) structured mesh
with $12$ by $12$ nodes (coarse) or $24$ by $24$ nodes (fine). The time step
varies between $0.2$ (coarse) and $0.1$ (fine).

\begin{figure}[h!]
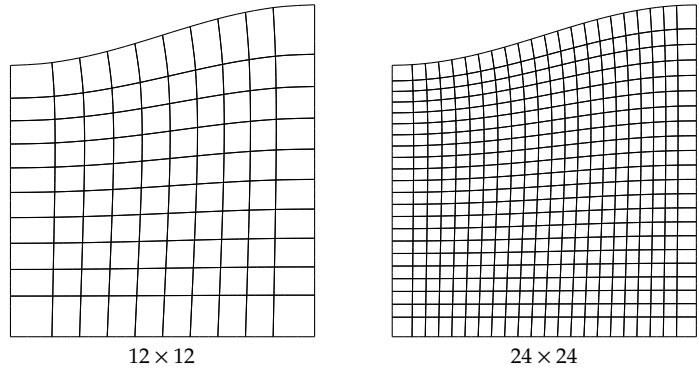

  \centering
  \begin{minipage}[b]{0.3\textwidth}
    \centering
    \begin{tikzpicture}[scale=4.0]
      \tikzstyle{mesh} = [thin]
      \input{Meshes/cos.iga.12x12.pdenor.dt100e-3/mesh0.tex}
    \end{tikzpicture}

    \footnotesize
    $12 \times 12$
  \end{minipage}
  \begin{minipage}[b]{0.3\textwidth}
    \centering
    \begin{tikzpicture}[scale=4.0]
      \tikzstyle{mesh} = [thin]
      \input{Meshes/cos.iga.24x24.pdenor.dt100e-3/mesh0.tex}
    \end{tikzpicture}

    \footnotesize
    $24 \times 24$
  \end{minipage}
  \caption{Meshes used for IGA simulations with $p=2$. Elements are shown in physical
           space.}
  \label{fig:tank.initialmeshes}
\end{figure}

The IGA simulations are conducted with a coarse mesh containing $12$ by $12$
control points or a fine mesh containing $24$ by $24$ control points. In both
cases, we consider a spline of degree 2. The element distribution in physical
space is indicated in Figure~\ref{fig:tank.initialmeshes}.

For both ---- standard FEM and IGA --- the initial position of the nodes /
control points describing the free surface is computed by solving a
least-squares fit problem.

\begin{figure}[h!]
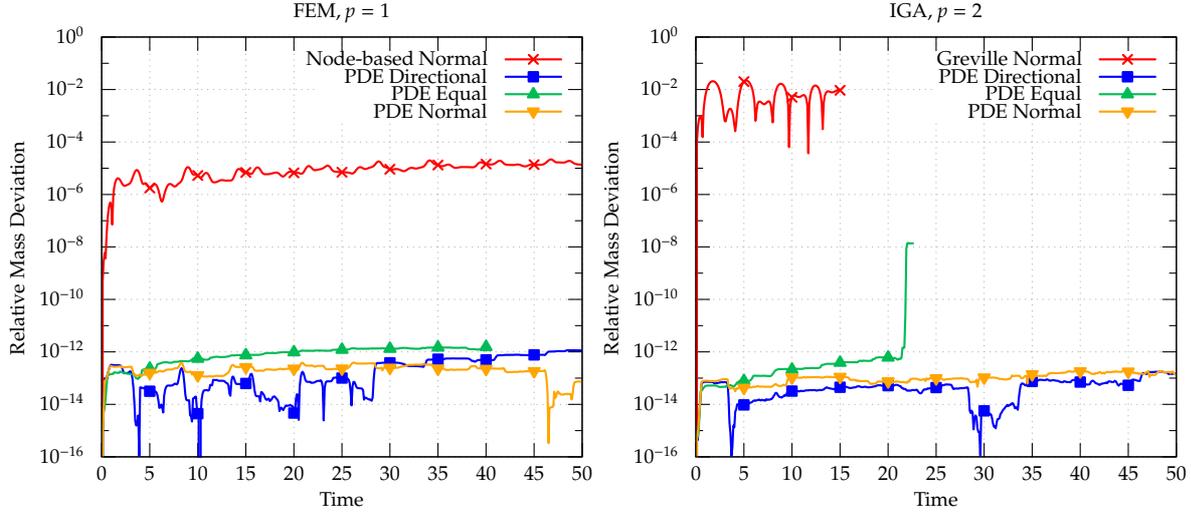

  \begin{minipage}[b]{0.47\textwidth}
    \centering
    \input{Plots/SloshingTank/mass.log.24x24.dt100e-3.fem.tex}
  \end{minipage}
  \begin{minipage}[b]{0.47\textwidth}
    \centering
    \input{Plots/SloshingTank/mass.log.24x24.dt100e-3.iga.tex}
  \end{minipage}
  \caption{Mass conservation error $|m(t)/m_0-1|$ for different displacement
methods of the free surface. On the left, the results for the standard finite
element method are collected. The new, PDE-based approaches show a significant
improvement of the mass conservation as compared to the traditional, node-based
versions. On the right, the results for the IGA versions are displayed. Again,
the PDE-versions preserve mass almost up to machine precision.}
  \label{fig:tank.mass}
\end{figure}

One standard means of comparison for this type of benchmark is the mass
conservation. In general, when utilizing a finite element method, mass loss is
expected. Figure~\ref{fig:tank.mass} plots the mass conservation error
${\tilde{m}}$, defined as

\begin{equation}
  {\tilde{m}}(t) \; := \; \left| \frac{m(t)}{m_0} - 1 \right| \; = \; \left| \frac{m(t)}{m(0)} - 1 \right| ,
\end{equation}

on a logarithmic scale. Here, $m(t)$ signifies the fluid mass in the tank, with
$m_0$ as initial mass. The left plot shows results for the standard finite
element method. Compared are (1) the traditional, node-based displacement in
normal direction and (2) the PDE-based imposition in the variants full velocity,
normal velocity and vertical velocity. The traditional version shows a mass loss
that noticeably deviates from the ideal value of $0.0$. Note that --- via design
--- the PDE-based imposition of the no-penetration boundary condition preserves
mass up to machine precision.

This observation still holds true for the IGA case displayed on the right plot.
Again, the point-based variant --- here, Greville points --- leads to a mass
loss, which is roughly nine orders of magnitude higher than the PDE-based
versions. Note that for both plots, a comparison of the PDE methods among each
other is not justifiable given how close the values are to machine precision.

In both plots, the simulations utilizing a displacement with the full fluid
velocity break off, before they reach the final simulation time. This is due to
the mesh becoming invalid (mesh tangling). This is again a substantiation of the
claim that movement in a fully Lagrangian sense is in general not useful. In the
IGA case, invalid meshes occur also for the Greville case.

\begin{figure}[h!]
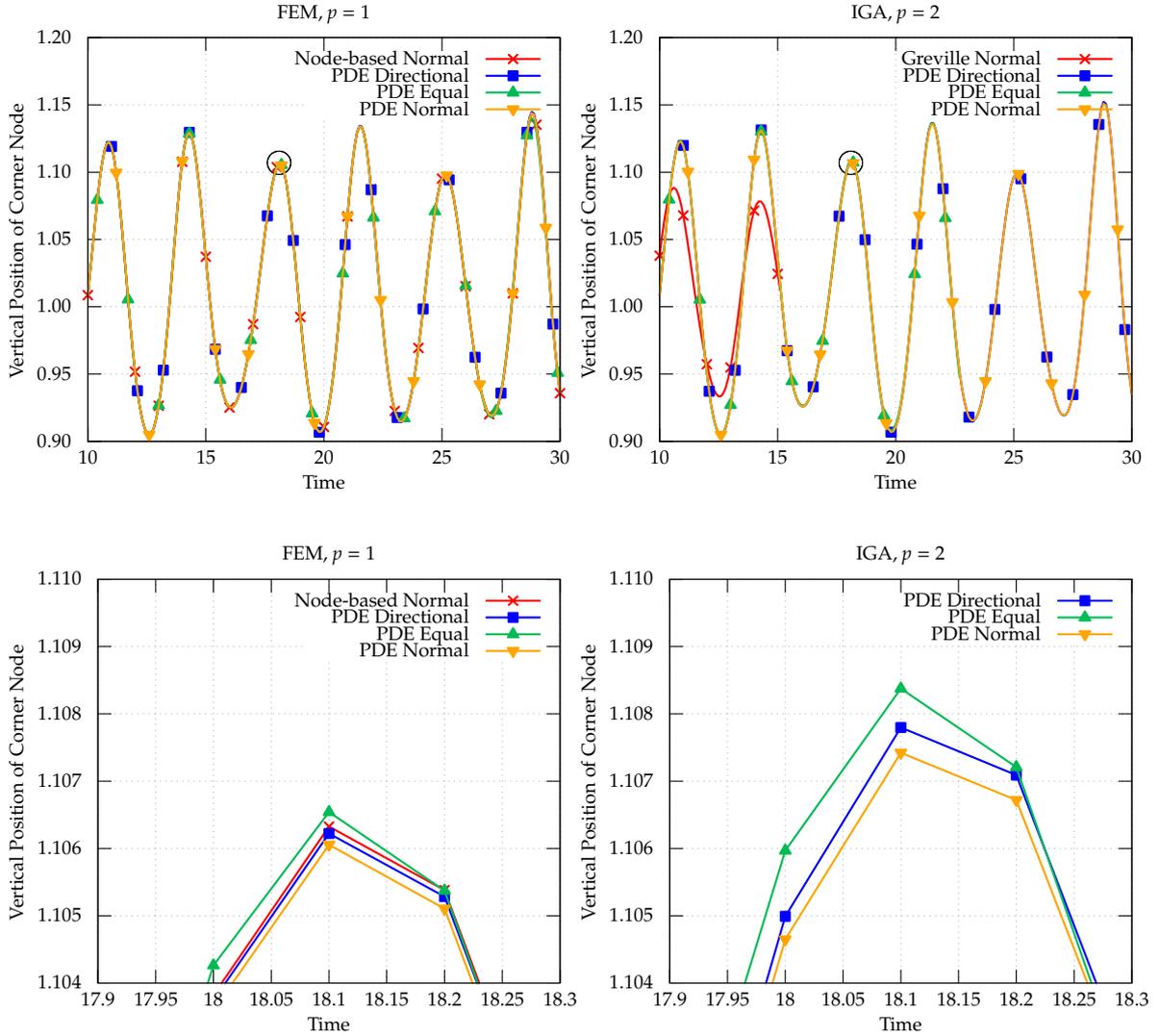

  \centering
  \begin{minipage}[b]{0.47\textwidth}
    \centering
    \input{Plots/SloshingTank/mov.24x24.dt100e-3.fem.tex}
  \end{minipage}
  \begin{minipage}[b]{0.47\textwidth}
    \centering
    \input{Plots/SloshingTank/mov.24x24.dt100e-3.iga.tex}
  \end{minipage}
  \begin{minipage}[b]{0.47\textwidth}
    \centering
    \begin{tikzpicture}[gnuplot]
\tikzset{every node/.append style={scale=0.75}}
\path (0.000,0.000) rectangle (8.000,7.000);
\gpcolor{color=gp lt color axes}
\gpsetlinetype{gp lt axes}
\gpsetdashtype{gp dt axes}
\gpsetlinewidth{0.50}
\draw[gp path] (1.266,0.739)--(7.585,0.739);
\gpcolor{color=gp lt color border}
\gpsetlinetype{gp lt border}
\gpsetdashtype{gp dt solid}
\gpsetlinewidth{1.00}
\draw[gp path] (1.266,0.739)--(1.446,0.739);
\draw[gp path] (7.585,0.739)--(7.405,0.739);
\node[gp node right] at (1.128,0.739) {1.104};
\gpcolor{color=gp lt color axes}
\gpsetlinetype{gp lt axes}
\gpsetdashtype{gp dt axes}
\gpsetlinewidth{0.50}
\draw[gp path] (1.266,1.667)--(7.585,1.667);
\gpcolor{color=gp lt color border}
\gpsetlinetype{gp lt border}
\gpsetdashtype{gp dt solid}
\gpsetlinewidth{1.00}
\draw[gp path] (1.266,1.667)--(1.446,1.667);
\draw[gp path] (7.585,1.667)--(7.405,1.667);
\node[gp node right] at (1.128,1.667) {1.105};
\gpcolor{color=gp lt color axes}
\gpsetlinetype{gp lt axes}
\gpsetdashtype{gp dt axes}
\gpsetlinewidth{0.50}
\draw[gp path] (1.266,2.595)--(7.585,2.595);
\gpcolor{color=gp lt color border}
\gpsetlinetype{gp lt border}
\gpsetdashtype{gp dt solid}
\gpsetlinewidth{1.00}
\draw[gp path] (1.266,2.595)--(1.446,2.595);
\draw[gp path] (7.585,2.595)--(7.405,2.595);
\node[gp node right] at (1.128,2.595) {1.106};
\gpcolor{color=gp lt color axes}
\gpsetlinetype{gp lt axes}
\gpsetdashtype{gp dt axes}
\gpsetlinewidth{0.50}
\draw[gp path] (1.266,3.522)--(7.585,3.522);
\gpcolor{color=gp lt color border}
\gpsetlinetype{gp lt border}
\gpsetdashtype{gp dt solid}
\gpsetlinewidth{1.00}
\draw[gp path] (1.266,3.522)--(1.446,3.522);
\draw[gp path] (7.585,3.522)--(7.405,3.522);
\node[gp node right] at (1.128,3.522) {1.107};
\gpcolor{color=gp lt color axes}
\gpsetlinetype{gp lt axes}
\gpsetdashtype{gp dt axes}
\gpsetlinewidth{0.50}
\draw[gp path] (1.266,4.450)--(7.585,4.450);
\gpcolor{color=gp lt color border}
\gpsetlinetype{gp lt border}
\gpsetdashtype{gp dt solid}
\gpsetlinewidth{1.00}
\draw[gp path] (1.266,4.450)--(1.446,4.450);
\draw[gp path] (7.585,4.450)--(7.405,4.450);
\node[gp node right] at (1.128,4.450) {1.108};
\gpcolor{color=gp lt color axes}
\gpsetlinetype{gp lt axes}
\gpsetdashtype{gp dt axes}
\gpsetlinewidth{0.50}
\draw[gp path] (1.266,5.378)--(4.093,5.378);
\draw[gp path] (7.447,5.378)--(7.585,5.378);
\gpcolor{color=gp lt color border}
\gpsetlinetype{gp lt border}
\gpsetdashtype{gp dt solid}
\gpsetlinewidth{1.00}
\draw[gp path] (1.266,5.378)--(1.446,5.378);
\draw[gp path] (7.585,5.378)--(7.405,5.378);
\node[gp node right] at (1.128,5.378) {1.109};
\gpcolor{color=gp lt color axes}
\gpsetlinetype{gp lt axes}
\gpsetdashtype{gp dt axes}
\gpsetlinewidth{0.50}
\draw[gp path] (1.266,6.306)--(7.585,6.306);
\gpcolor{color=gp lt color border}
\gpsetlinetype{gp lt border}
\gpsetdashtype{gp dt solid}
\gpsetlinewidth{1.00}
\draw[gp path] (1.266,6.306)--(1.446,6.306);
\draw[gp path] (7.585,6.306)--(7.405,6.306);
\node[gp node right] at (1.128,6.306) {1.110};
\gpcolor{color=gp lt color axes}
\gpsetlinetype{gp lt axes}
\gpsetdashtype{gp dt axes}
\gpsetlinewidth{0.50}
\draw[gp path] (1.266,0.739)--(1.266,6.306);
\gpcolor{color=gp lt color border}
\gpsetlinetype{gp lt border}
\gpsetdashtype{gp dt solid}
\gpsetlinewidth{1.00}
\draw[gp path] (1.266,0.739)--(1.266,0.919);
\draw[gp path] (1.266,6.306)--(1.266,6.126);
\node[gp node center] at (1.266,0.508) {$17.9$};
\gpcolor{color=gp lt color axes}
\gpsetlinetype{gp lt axes}
\gpsetdashtype{gp dt axes}
\gpsetlinewidth{0.50}
\draw[gp path] (2.056,0.739)--(2.056,6.306);
\gpcolor{color=gp lt color border}
\gpsetlinetype{gp lt border}
\gpsetdashtype{gp dt solid}
\gpsetlinewidth{1.00}
\draw[gp path] (2.056,0.739)--(2.056,0.919);
\draw[gp path] (2.056,6.306)--(2.056,6.126);
\node[gp node center] at (2.056,0.508) {$17.95$};
\gpcolor{color=gp lt color axes}
\gpsetlinetype{gp lt axes}
\gpsetdashtype{gp dt axes}
\gpsetlinewidth{0.50}
\draw[gp path] (2.846,0.739)--(2.846,6.306);
\gpcolor{color=gp lt color border}
\gpsetlinetype{gp lt border}
\gpsetdashtype{gp dt solid}
\gpsetlinewidth{1.00}
\draw[gp path] (2.846,0.739)--(2.846,0.919);
\draw[gp path] (2.846,6.306)--(2.846,6.126);
\node[gp node center] at (2.846,0.508) {$18$};
\gpcolor{color=gp lt color axes}
\gpsetlinetype{gp lt axes}
\gpsetdashtype{gp dt axes}
\gpsetlinewidth{0.50}
\draw[gp path] (3.636,0.739)--(3.636,6.306);
\gpcolor{color=gp lt color border}
\gpsetlinetype{gp lt border}
\gpsetdashtype{gp dt solid}
\gpsetlinewidth{1.00}
\draw[gp path] (3.636,0.739)--(3.636,0.919);
\draw[gp path] (3.636,6.306)--(3.636,6.126);
\node[gp node center] at (3.636,0.508) {$18.05$};
\gpcolor{color=gp lt color axes}
\gpsetlinetype{gp lt axes}
\gpsetdashtype{gp dt axes}
\gpsetlinewidth{0.50}
\draw[gp path] (4.426,0.739)--(4.426,5.202);
\draw[gp path] (4.426,6.126)--(4.426,6.306);
\gpcolor{color=gp lt color border}
\gpsetlinetype{gp lt border}
\gpsetdashtype{gp dt solid}
\gpsetlinewidth{1.00}
\draw[gp path] (4.426,0.739)--(4.426,0.919);
\draw[gp path] (4.426,6.306)--(4.426,6.126);
\node[gp node center] at (4.426,0.508) {$18.1$};
\gpcolor{color=gp lt color axes}
\gpsetlinetype{gp lt axes}
\gpsetdashtype{gp dt axes}
\gpsetlinewidth{0.50}
\draw[gp path] (5.215,0.739)--(5.215,5.202);
\draw[gp path] (5.215,6.126)--(5.215,6.306);
\gpcolor{color=gp lt color border}
\gpsetlinetype{gp lt border}
\gpsetdashtype{gp dt solid}
\gpsetlinewidth{1.00}
\draw[gp path] (5.215,0.739)--(5.215,0.919);
\draw[gp path] (5.215,6.306)--(5.215,6.126);
\node[gp node center] at (5.215,0.508) {$18.15$};
\gpcolor{color=gp lt color axes}
\gpsetlinetype{gp lt axes}
\gpsetdashtype{gp dt axes}
\gpsetlinewidth{0.50}
\draw[gp path] (6.005,0.739)--(6.005,5.202);
\draw[gp path] (6.005,6.126)--(6.005,6.306);
\gpcolor{color=gp lt color border}
\gpsetlinetype{gp lt border}
\gpsetdashtype{gp dt solid}
\gpsetlinewidth{1.00}
\draw[gp path] (6.005,0.739)--(6.005,0.919);
\draw[gp path] (6.005,6.306)--(6.005,6.126);
\node[gp node center] at (6.005,0.508) {$18.2$};
\gpcolor{color=gp lt color axes}
\gpsetlinetype{gp lt axes}
\gpsetdashtype{gp dt axes}
\gpsetlinewidth{0.50}
\draw[gp path] (6.795,0.739)--(6.795,5.202);
\draw[gp path] (6.795,6.126)--(6.795,6.306);
\gpcolor{color=gp lt color border}
\gpsetlinetype{gp lt border}
\gpsetdashtype{gp dt solid}
\gpsetlinewidth{1.00}
\draw[gp path] (6.795,0.739)--(6.795,0.919);
\draw[gp path] (6.795,6.306)--(6.795,6.126);
\node[gp node center] at (6.795,0.508) {$18.25$};
\gpcolor{color=gp lt color axes}
\gpsetlinetype{gp lt axes}
\gpsetdashtype{gp dt axes}
\gpsetlinewidth{0.50}
\draw[gp path] (7.585,0.739)--(7.585,6.306);
\gpcolor{color=gp lt color border}
\gpsetlinetype{gp lt border}
\gpsetdashtype{gp dt solid}
\gpsetlinewidth{1.00}
\draw[gp path] (7.585,0.739)--(7.585,0.919);
\draw[gp path] (7.585,6.306)--(7.585,6.126);
\node[gp node center] at (7.585,0.508) {$18.3$};
\draw[gp path] (1.266,6.306)--(1.266,0.739)--(7.585,0.739)--(7.585,6.306)--cycle;
\node[gp node center,rotate=-270] at (0.184,3.522) {Vertical Position of Corner Node};
\node[gp node center] at (4.425,0.162) {Time};
\node[gp node center] at (4.425,6.653) {FEM, $p=1$};
\node[gp node right] at (6.439,6.010) {Node-based Normal};
\gpcolor{rgb color={1.000,0.000,0.000}}
\gpsetlinewidth{2.00}
\draw[gp path] (6.577,6.010)--(7.309,6.010);
\draw[gp path] (2.958,0.739)--(4.426,2.897)--(6.005,2.021)--(6.480,0.739);
\gpsetpointsize{6.00}
\gppoint{gp mark 2}{(4.426,2.897)}
\gppoint{gp mark 2}{(6.005,2.021)}
\gppoint{gp mark 2}{(6.943,6.010)}
\gpcolor{color=gp lt color border}
\node[gp node right] at (6.439,5.779) {PDE Directional};
\gpcolor{rgb color={0.000,0.000,1.000}}
\draw[gp path] (6.577,5.779)--(7.309,5.779);
\draw[gp path] (3.008,0.739)--(4.426,2.801)--(6.005,1.931)--(6.451,0.739);
\gpsetpointsize{4.40}
\gppoint{gp mark 5}{(4.426,2.801)}
\gppoint{gp mark 5}{(6.005,1.931)}
\gppoint{gp mark 5}{(6.943,5.779)}
\gpcolor{color=gp lt color border}
\node[gp node right] at (6.439,5.548) {PDE Equal};
\gpcolor{rgb color={0.000,0.733,0.333}}
\draw[gp path] (6.577,5.548)--(7.309,5.548);
\draw[gp path] (2.769,0.739)--(2.846,0.983)--(4.426,3.097)--(6.005,2.014)--(6.458,0.739);
\gpsetpointsize{6.00}
\gppoint{gp mark 9}{(2.846,0.983)}
\gppoint{gp mark 9}{(4.426,3.097)}
\gppoint{gp mark 9}{(6.005,2.014)}
\gppoint{gp mark 9}{(6.943,5.548)}
\gpcolor{color=gp lt color border}
\node[gp node right] at (6.439,5.317) {PDE Normal};
\gpcolor{rgb color={1.000,0.647,0.000}}
\draw[gp path] (6.577,5.317)--(7.309,5.317);
\draw[gp path] (3.106,0.739)--(4.426,2.646)--(6.005,1.765)--(6.389,0.739);
\gppoint{gp mark 11}{(4.426,2.646)}
\gppoint{gp mark 11}{(6.005,1.765)}
\gppoint{gp mark 11}{(6.943,5.317)}
\gpcolor{color=gp lt color border}
\gpsetlinewidth{1.00}
\draw[gp path] (1.266,6.306)--(1.266,0.739)--(7.585,0.739)--(7.585,6.306)--cycle;
\gpdefrectangularnode{gp plot 1}{\pgfpoint{1.266cm}{0.739cm}}{\pgfpoint{7.585cm}{6.306cm}}
\end{tikzpicture}
  \end{minipage}
  \begin{minipage}[b]{0.47\textwidth}
    \centering
    \begin{tikzpicture}[gnuplot]
\tikzset{every node/.append style={scale=0.75}}
\path (0.000,0.000) rectangle (8.000,7.000);
\gpcolor{color=gp lt color axes}
\gpsetlinetype{gp lt axes}
\gpsetdashtype{gp dt axes}
\gpsetlinewidth{0.50}
\draw[gp path] (1.266,0.739)--(7.585,0.739);
\gpcolor{color=gp lt color border}
\gpsetlinetype{gp lt border}
\gpsetdashtype{gp dt solid}
\gpsetlinewidth{1.00}
\draw[gp path] (1.266,0.739)--(1.446,0.739);
\draw[gp path] (7.585,0.739)--(7.405,0.739);
\node[gp node right] at (1.128,0.739) {1.104};
\gpcolor{color=gp lt color axes}
\gpsetlinetype{gp lt axes}
\gpsetdashtype{gp dt axes}
\gpsetlinewidth{0.50}
\draw[gp path] (1.266,1.667)--(7.585,1.667);
\gpcolor{color=gp lt color border}
\gpsetlinetype{gp lt border}
\gpsetdashtype{gp dt solid}
\gpsetlinewidth{1.00}
\draw[gp path] (1.266,1.667)--(1.446,1.667);
\draw[gp path] (7.585,1.667)--(7.405,1.667);
\node[gp node right] at (1.128,1.667) {1.105};
\gpcolor{color=gp lt color axes}
\gpsetlinetype{gp lt axes}
\gpsetdashtype{gp dt axes}
\gpsetlinewidth{0.50}
\draw[gp path] (1.266,2.595)--(7.585,2.595);
\gpcolor{color=gp lt color border}
\gpsetlinetype{gp lt border}
\gpsetdashtype{gp dt solid}
\gpsetlinewidth{1.00}
\draw[gp path] (1.266,2.595)--(1.446,2.595);
\draw[gp path] (7.585,2.595)--(7.405,2.595);
\node[gp node right] at (1.128,2.595) {1.106};
\gpcolor{color=gp lt color axes}
\gpsetlinetype{gp lt axes}
\gpsetdashtype{gp dt axes}
\gpsetlinewidth{0.50}
\draw[gp path] (1.266,3.522)--(7.585,3.522);
\gpcolor{color=gp lt color border}
\gpsetlinetype{gp lt border}
\gpsetdashtype{gp dt solid}
\gpsetlinewidth{1.00}
\draw[gp path] (1.266,3.522)--(1.446,3.522);
\draw[gp path] (7.585,3.522)--(7.405,3.522);
\node[gp node right] at (1.128,3.522) {1.107};
\gpcolor{color=gp lt color axes}
\gpsetlinetype{gp lt axes}
\gpsetdashtype{gp dt axes}
\gpsetlinewidth{0.50}
\draw[gp path] (1.266,4.450)--(7.585,4.450);
\gpcolor{color=gp lt color border}
\gpsetlinetype{gp lt border}
\gpsetdashtype{gp dt solid}
\gpsetlinewidth{1.00}
\draw[gp path] (1.266,4.450)--(1.446,4.450);
\draw[gp path] (7.585,4.450)--(7.405,4.450);
\node[gp node right] at (1.128,4.450) {1.108};
\gpcolor{color=gp lt color axes}
\gpsetlinetype{gp lt axes}
\gpsetdashtype{gp dt axes}
\gpsetlinewidth{0.50}
\draw[gp path] (1.266,5.378)--(7.585,5.378);
\gpcolor{color=gp lt color border}
\gpsetlinetype{gp lt border}
\gpsetdashtype{gp dt solid}
\gpsetlinewidth{1.00}
\draw[gp path] (1.266,5.378)--(1.446,5.378);
\draw[gp path] (7.585,5.378)--(7.405,5.378);
\node[gp node right] at (1.128,5.378) {1.109};
\gpcolor{color=gp lt color axes}
\gpsetlinetype{gp lt axes}
\gpsetdashtype{gp dt axes}
\gpsetlinewidth{0.50}
\draw[gp path] (1.266,6.306)--(7.585,6.306);
\gpcolor{color=gp lt color border}
\gpsetlinetype{gp lt border}
\gpsetdashtype{gp dt solid}
\gpsetlinewidth{1.00}
\draw[gp path] (1.266,6.306)--(1.446,6.306);
\draw[gp path] (7.585,6.306)--(7.405,6.306);
\node[gp node right] at (1.128,6.306) {1.110};
\gpcolor{color=gp lt color axes}
\gpsetlinetype{gp lt axes}
\gpsetdashtype{gp dt axes}
\gpsetlinewidth{0.50}
\draw[gp path] (1.266,0.739)--(1.266,6.306);
\gpcolor{color=gp lt color border}
\gpsetlinetype{gp lt border}
\gpsetdashtype{gp dt solid}
\gpsetlinewidth{1.00}
\draw[gp path] (1.266,0.739)--(1.266,0.919);
\draw[gp path] (1.266,6.306)--(1.266,6.126);
\node[gp node center] at (1.266,0.508) {$17.9$};
\gpcolor{color=gp lt color axes}
\gpsetlinetype{gp lt axes}
\gpsetdashtype{gp dt axes}
\gpsetlinewidth{0.50}
\draw[gp path] (2.056,0.739)--(2.056,6.306);
\gpcolor{color=gp lt color border}
\gpsetlinetype{gp lt border}
\gpsetdashtype{gp dt solid}
\gpsetlinewidth{1.00}
\draw[gp path] (2.056,0.739)--(2.056,0.919);
\draw[gp path] (2.056,6.306)--(2.056,6.126);
\node[gp node center] at (2.056,0.508) {$17.95$};
\gpcolor{color=gp lt color axes}
\gpsetlinetype{gp lt axes}
\gpsetdashtype{gp dt axes}
\gpsetlinewidth{0.50}
\draw[gp path] (2.846,0.739)--(2.846,6.306);
\gpcolor{color=gp lt color border}
\gpsetlinetype{gp lt border}
\gpsetdashtype{gp dt solid}
\gpsetlinewidth{1.00}
\draw[gp path] (2.846,0.739)--(2.846,0.919);
\draw[gp path] (2.846,6.306)--(2.846,6.126);
\node[gp node center] at (2.846,0.508) {$18$};
\gpcolor{color=gp lt color axes}
\gpsetlinetype{gp lt axes}
\gpsetdashtype{gp dt axes}
\gpsetlinewidth{0.50}
\draw[gp path] (3.636,0.739)--(3.636,6.306);
\gpcolor{color=gp lt color border}
\gpsetlinetype{gp lt border}
\gpsetdashtype{gp dt solid}
\gpsetlinewidth{1.00}
\draw[gp path] (3.636,0.739)--(3.636,0.919);
\draw[gp path] (3.636,6.306)--(3.636,6.126);
\node[gp node center] at (3.636,0.508) {$18.05$};
\gpcolor{color=gp lt color axes}
\gpsetlinetype{gp lt axes}
\gpsetdashtype{gp dt axes}
\gpsetlinewidth{0.50}
\draw[gp path] (4.426,0.739)--(4.426,5.433);
\draw[gp path] (4.426,6.126)--(4.426,6.306);
\gpcolor{color=gp lt color border}
\gpsetlinetype{gp lt border}
\gpsetdashtype{gp dt solid}
\gpsetlinewidth{1.00}
\draw[gp path] (4.426,0.739)--(4.426,0.919);
\draw[gp path] (4.426,6.306)--(4.426,6.126);
\node[gp node center] at (4.426,0.508) {$18.1$};
\gpcolor{color=gp lt color axes}
\gpsetlinetype{gp lt axes}
\gpsetdashtype{gp dt axes}
\gpsetlinewidth{0.50}
\draw[gp path] (5.215,0.739)--(5.215,5.433);
\draw[gp path] (5.215,6.126)--(5.215,6.306);
\gpcolor{color=gp lt color border}
\gpsetlinetype{gp lt border}
\gpsetdashtype{gp dt solid}
\gpsetlinewidth{1.00}
\draw[gp path] (5.215,0.739)--(5.215,0.919);
\draw[gp path] (5.215,6.306)--(5.215,6.126);
\node[gp node center] at (5.215,0.508) {$18.15$};
\gpcolor{color=gp lt color axes}
\gpsetlinetype{gp lt axes}
\gpsetdashtype{gp dt axes}
\gpsetlinewidth{0.50}
\draw[gp path] (6.005,0.739)--(6.005,5.433);
\draw[gp path] (6.005,6.126)--(6.005,6.306);
\gpcolor{color=gp lt color border}
\gpsetlinetype{gp lt border}
\gpsetdashtype{gp dt solid}
\gpsetlinewidth{1.00}
\draw[gp path] (6.005,0.739)--(6.005,0.919);
\draw[gp path] (6.005,6.306)--(6.005,6.126);
\node[gp node center] at (6.005,0.508) {$18.2$};
\gpcolor{color=gp lt color axes}
\gpsetlinetype{gp lt axes}
\gpsetdashtype{gp dt axes}
\gpsetlinewidth{0.50}
\draw[gp path] (6.795,0.739)--(6.795,5.433);
\draw[gp path] (6.795,6.126)--(6.795,6.306);
\gpcolor{color=gp lt color border}
\gpsetlinetype{gp lt border}
\gpsetdashtype{gp dt solid}
\gpsetlinewidth{1.00}
\draw[gp path] (6.795,0.739)--(6.795,0.919);
\draw[gp path] (6.795,6.306)--(6.795,6.126);
\node[gp node center] at (6.795,0.508) {$18.25$};
\gpcolor{color=gp lt color axes}
\gpsetlinetype{gp lt axes}
\gpsetdashtype{gp dt axes}
\gpsetlinewidth{0.50}
\draw[gp path] (7.585,0.739)--(7.585,6.306);
\gpcolor{color=gp lt color border}
\gpsetlinetype{gp lt border}
\gpsetdashtype{gp dt solid}
\gpsetlinewidth{1.00}
\draw[gp path] (7.585,0.739)--(7.585,0.919);
\draw[gp path] (7.585,6.306)--(7.585,6.126);
\node[gp node center] at (7.585,0.508) {$18.3$};
\draw[gp path] (1.266,6.306)--(1.266,0.739)--(7.585,0.739)--(7.585,6.306)--cycle;
\node[gp node center,rotate=-270] at (0.184,3.522) {Vertical Position of Corner Node};
\node[gp node center] at (4.425,0.162) {Time};
\node[gp node center] at (4.425,6.653) {IGA, $p=2$};
\node[gp node right] at (6.439,6.010) {PDE Directional};
\gpcolor{rgb color={0.000,0.000,1.000}}
\gpsetlinewidth{2.00}
\draw[gp path] (6.577,6.010)--(7.309,6.010);
\draw[gp path] (2.582,0.739)--(2.846,1.659)--(4.426,4.262)--(6.005,3.606)--(7.098,0.739);
\gpsetpointsize{4.40}
\gppoint{gp mark 5}{(2.846,1.659)}
\gppoint{gp mark 5}{(4.426,4.262)}
\gppoint{gp mark 5}{(6.005,3.606)}
\gppoint{gp mark 5}{(6.943,6.010)}
\gpcolor{color=gp lt color border}
\node[gp node right] at (6.439,5.779) {PDE Equal};
\gpcolor{rgb color={0.000,0.733,0.333}}
\draw[gp path] (6.577,5.779)--(7.309,5.779);
\draw[gp path] (2.294,0.739)--(2.846,2.566)--(4.426,4.799)--(6.005,3.716)--(7.032,0.739);
\gpsetpointsize{6.00}
\gppoint{gp mark 9}{(2.846,2.566)}
\gppoint{gp mark 9}{(4.426,4.799)}
\gppoint{gp mark 9}{(6.005,3.716)}
\gppoint{gp mark 9}{(6.943,5.779)}
\gpcolor{color=gp lt color border}
\node[gp node right] at (6.439,5.548) {PDE Normal};
\gpcolor{rgb color={1.000,0.647,0.000}}
\draw[gp path] (6.577,5.548)--(7.309,5.548);
\draw[gp path] (2.671,0.739)--(2.846,1.342)--(4.426,3.914)--(6.005,3.262)--(6.977,0.739);
\gppoint{gp mark 11}{(2.846,1.342)}
\gppoint{gp mark 11}{(4.426,3.914)}
\gppoint{gp mark 11}{(6.005,3.262)}
\gppoint{gp mark 11}{(6.943,5.548)}
\gpcolor{color=gp lt color border}
\gpsetlinewidth{1.00}
\draw[gp path] (1.266,6.306)--(1.266,0.739)--(7.585,0.739)--(7.585,6.306)--cycle;
\gpdefrectangularnode{gp plot 1}{\pgfpoint{1.266cm}{0.739cm}}{\pgfpoint{7.585cm}{6.306cm}}
\end{tikzpicture}
  \end{minipage}
  \caption{As an indicator of the free-surface shape, the vertical position of
upper right corner of the computational domain is plotted over time for linear
FEM (left) and quadratic IGA (right). The lower plots constitute a zoom on the
respective upper plot. All simulations were conducted on the fine mesh.}
  \label{fig:tank.mov}
\end{figure}

\begin{figure}[h!]
  \centering
  \begin{minipage}[b]{0.47\textwidth}
    \begin{tikzpicture}[gnuplot]
\tikzset{every node/.append style={scale=0.75}}
\path (0.000,0.000) rectangle (8.000,7.000);
\gpcolor{color=gp lt color axes}
\gpsetlinetype{gp lt axes}
\gpsetdashtype{gp dt axes}
\gpsetlinewidth{0.50}
\draw[gp path] (1.128,0.739)--(7.585,0.739);
\gpcolor{color=gp lt color border}
\gpsetlinetype{gp lt border}
\gpsetdashtype{gp dt solid}
\gpsetlinewidth{1.00}
\draw[gp path] (1.128,0.739)--(1.308,0.739);
\draw[gp path] (7.585,0.739)--(7.405,0.739);
\node[gp node right] at (0.990,0.739) {0.90};
\gpcolor{color=gp lt color axes}
\gpsetlinetype{gp lt axes}
\gpsetdashtype{gp dt axes}
\gpsetlinewidth{0.50}
\draw[gp path] (1.128,1.667)--(7.585,1.667);
\gpcolor{color=gp lt color border}
\gpsetlinetype{gp lt border}
\gpsetdashtype{gp dt solid}
\gpsetlinewidth{1.00}
\draw[gp path] (1.128,1.667)--(1.308,1.667);
\draw[gp path] (7.585,1.667)--(7.405,1.667);
\node[gp node right] at (0.990,1.667) {0.95};
\gpcolor{color=gp lt color axes}
\gpsetlinetype{gp lt axes}
\gpsetdashtype{gp dt axes}
\gpsetlinewidth{0.50}
\draw[gp path] (1.128,2.595)--(7.585,2.595);
\gpcolor{color=gp lt color border}
\gpsetlinetype{gp lt border}
\gpsetdashtype{gp dt solid}
\gpsetlinewidth{1.00}
\draw[gp path] (1.128,2.595)--(1.308,2.595);
\draw[gp path] (7.585,2.595)--(7.405,2.595);
\node[gp node right] at (0.990,2.595) {1.00};
\gpcolor{color=gp lt color axes}
\gpsetlinetype{gp lt axes}
\gpsetdashtype{gp dt axes}
\gpsetlinewidth{0.50}
\draw[gp path] (1.128,3.523)--(7.585,3.523);
\gpcolor{color=gp lt color border}
\gpsetlinetype{gp lt border}
\gpsetdashtype{gp dt solid}
\gpsetlinewidth{1.00}
\draw[gp path] (1.128,3.523)--(1.308,3.523);
\draw[gp path] (7.585,3.523)--(7.405,3.523);
\node[gp node right] at (0.990,3.523) {1.05};
\gpcolor{color=gp lt color axes}
\gpsetlinetype{gp lt axes}
\gpsetdashtype{gp dt axes}
\gpsetlinewidth{0.50}
\draw[gp path] (1.128,4.450)--(7.585,4.450);
\gpcolor{color=gp lt color border}
\gpsetlinetype{gp lt border}
\gpsetdashtype{gp dt solid}
\gpsetlinewidth{1.00}
\draw[gp path] (1.128,4.450)--(1.308,4.450);
\draw[gp path] (7.585,4.450)--(7.405,4.450);
\node[gp node right] at (0.990,4.450) {1.10};
\gpcolor{color=gp lt color axes}
\gpsetlinetype{gp lt axes}
\gpsetdashtype{gp dt axes}
\gpsetlinewidth{0.50}
\draw[gp path] (1.128,5.378)--(7.585,5.378);
\gpcolor{color=gp lt color border}
\gpsetlinetype{gp lt border}
\gpsetdashtype{gp dt solid}
\gpsetlinewidth{1.00}
\draw[gp path] (1.128,5.378)--(1.308,5.378);
\draw[gp path] (7.585,5.378)--(7.405,5.378);
\node[gp node right] at (0.990,5.378) {1.15};
\gpcolor{color=gp lt color axes}
\gpsetlinetype{gp lt axes}
\gpsetdashtype{gp dt axes}
\gpsetlinewidth{0.50}
\draw[gp path] (1.128,6.306)--(7.585,6.306);
\gpcolor{color=gp lt color border}
\gpsetlinetype{gp lt border}
\gpsetdashtype{gp dt solid}
\gpsetlinewidth{1.00}
\draw[gp path] (1.128,6.306)--(1.308,6.306);
\draw[gp path] (7.585,6.306)--(7.405,6.306);
\node[gp node right] at (0.990,6.306) {1.20};
\gpcolor{color=gp lt color axes}
\gpsetlinetype{gp lt axes}
\gpsetdashtype{gp dt axes}
\gpsetlinewidth{0.50}
\draw[gp path] (1.128,0.739)--(1.128,6.306);
\gpcolor{color=gp lt color border}
\gpsetlinetype{gp lt border}
\gpsetdashtype{gp dt solid}
\gpsetlinewidth{1.00}
\draw[gp path] (1.128,0.739)--(1.128,0.919);
\draw[gp path] (1.128,6.306)--(1.128,6.126);
\node[gp node center] at (1.128,0.508) {$42$};
\gpcolor{color=gp lt color axes}
\gpsetlinetype{gp lt axes}
\gpsetdashtype{gp dt axes}
\gpsetlinewidth{0.50}
\draw[gp path] (2.050,0.739)--(2.050,6.306);
\gpcolor{color=gp lt color border}
\gpsetlinetype{gp lt border}
\gpsetdashtype{gp dt solid}
\gpsetlinewidth{1.00}
\draw[gp path] (2.050,0.739)--(2.050,0.919);
\draw[gp path] (2.050,6.306)--(2.050,6.126);
\node[gp node center] at (2.050,0.508) {$43$};
\gpcolor{color=gp lt color axes}
\gpsetlinetype{gp lt axes}
\gpsetdashtype{gp dt axes}
\gpsetlinewidth{0.50}
\draw[gp path] (2.973,0.739)--(2.973,6.306);
\gpcolor{color=gp lt color border}
\gpsetlinetype{gp lt border}
\gpsetdashtype{gp dt solid}
\gpsetlinewidth{1.00}
\draw[gp path] (2.973,0.739)--(2.973,0.919);
\draw[gp path] (2.973,6.306)--(2.973,6.126);
\node[gp node center] at (2.973,0.508) {$44$};
\gpcolor{color=gp lt color axes}
\gpsetlinetype{gp lt axes}
\gpsetdashtype{gp dt axes}
\gpsetlinewidth{0.50}
\draw[gp path] (3.895,0.739)--(3.895,6.306);
\gpcolor{color=gp lt color border}
\gpsetlinetype{gp lt border}
\gpsetdashtype{gp dt solid}
\gpsetlinewidth{1.00}
\draw[gp path] (3.895,0.739)--(3.895,0.919);
\draw[gp path] (3.895,6.306)--(3.895,6.126);
\node[gp node center] at (3.895,0.508) {$45$};
\gpcolor{color=gp lt color axes}
\gpsetlinetype{gp lt axes}
\gpsetdashtype{gp dt axes}
\gpsetlinewidth{0.50}
\draw[gp path] (4.818,0.739)--(4.818,6.306);
\gpcolor{color=gp lt color border}
\gpsetlinetype{gp lt border}
\gpsetdashtype{gp dt solid}
\gpsetlinewidth{1.00}
\draw[gp path] (4.818,0.739)--(4.818,0.919);
\draw[gp path] (4.818,6.306)--(4.818,6.126);
\node[gp node center] at (4.818,0.508) {$46$};
\gpcolor{color=gp lt color axes}
\gpsetlinetype{gp lt axes}
\gpsetdashtype{gp dt axes}
\gpsetlinewidth{0.50}
\draw[gp path] (5.740,0.739)--(5.740,5.664);
\draw[gp path] (5.740,6.126)--(5.740,6.306);
\gpcolor{color=gp lt color border}
\gpsetlinetype{gp lt border}
\gpsetdashtype{gp dt solid}
\gpsetlinewidth{1.00}
\draw[gp path] (5.740,0.739)--(5.740,0.919);
\draw[gp path] (5.740,6.306)--(5.740,6.126);
\node[gp node center] at (5.740,0.508) {$47$};
\gpcolor{color=gp lt color axes}
\gpsetlinetype{gp lt axes}
\gpsetdashtype{gp dt axes}
\gpsetlinewidth{0.50}
\draw[gp path] (6.663,0.739)--(6.663,5.664);
\draw[gp path] (6.663,6.126)--(6.663,6.306);
\gpcolor{color=gp lt color border}
\gpsetlinetype{gp lt border}
\gpsetdashtype{gp dt solid}
\gpsetlinewidth{1.00}
\draw[gp path] (6.663,0.739)--(6.663,0.919);
\draw[gp path] (6.663,6.306)--(6.663,6.126);
\node[gp node center] at (6.663,0.508) {$48$};
\gpcolor{color=gp lt color axes}
\gpsetlinetype{gp lt axes}
\gpsetdashtype{gp dt axes}
\gpsetlinewidth{0.50}
\draw[gp path] (7.585,0.739)--(7.585,6.306);
\gpcolor{color=gp lt color border}
\gpsetlinetype{gp lt border}
\gpsetdashtype{gp dt solid}
\gpsetlinewidth{1.00}
\draw[gp path] (7.585,0.739)--(7.585,0.919);
\draw[gp path] (7.585,6.306)--(7.585,6.126);
\node[gp node center] at (7.585,0.508) {$49$};
\draw[gp path] (1.128,6.306)--(1.128,0.739)--(7.585,0.739)--(7.585,6.306)--cycle;
\node[gp node center,rotate=-270] at (0.184,3.522) {Vertical Position of Corner Node};
\node[gp node center] at (4.356,0.162) {Time};
\node[gp node center] at (4.356,6.653) {$12 \times 12$ nodes, $\Delta t = 0.2$};
\node[gp node right] at (6.439,6.010) {FEM, $p=1$};
\gpcolor{rgb color={0.000,0.000,0.000}}
\gpsetlinewidth{2.00}
\draw[gp path] (6.577,6.010)--(7.309,6.010);
\draw[gp path] (1.128,1.577)--(1.312,1.974)--(1.497,2.457)--(1.681,2.988)--(1.866,3.513)%
  --(2.050,3.966)--(2.235,4.280)--(2.419,4.401)--(2.604,4.306)--(2.788,4.002)--(2.973,3.534)%
  --(3.157,2.972)--(3.342,2.397)--(3.526,1.877)--(3.711,1.464)--(3.895,1.186)--(4.080,1.062)%
  --(4.264,1.102)--(4.449,1.311)--(4.633,1.691)--(4.818,2.231)--(5.002,2.896)--(5.187,3.613)%
  --(5.371,4.270)--(5.556,4.734)--(5.740,4.897)--(5.925,4.717)--(6.109,4.232)--(6.294,3.552)%
  --(6.478,2.812)--(6.663,2.132)--(6.847,1.589)--(7.032,1.219)--(7.216,1.032)--(7.401,1.023)%
  --(7.585,1.186);
\gpsetpointsize{6.00}
\gppoint{gp mark 2}{(1.128,1.577)}
\gppoint{gp mark 2}{(1.681,2.988)}
\gppoint{gp mark 2}{(2.235,4.280)}
\gppoint{gp mark 2}{(2.788,4.002)}
\gppoint{gp mark 2}{(3.342,2.397)}
\gppoint{gp mark 2}{(3.895,1.186)}
\gppoint{gp mark 2}{(4.449,1.311)}
\gppoint{gp mark 2}{(5.002,2.896)}
\gppoint{gp mark 2}{(5.556,4.734)}
\gppoint{gp mark 2}{(6.109,4.232)}
\gppoint{gp mark 2}{(6.663,2.132)}
\gppoint{gp mark 2}{(7.216,1.032)}
\gppoint{gp mark 2}{(6.943,6.010)}
\gpcolor{color=gp lt color border}
\node[gp node right] at (6.439,5.779) {IGA, $p=2$};
\gpcolor{rgb color={0.000,0.000,1.000}}
\draw[gp path] (6.577,5.779)--(7.309,5.779);
\draw[gp path] (1.128,1.995)--(1.312,2.373)--(1.497,2.826)--(1.681,3.328)--(1.866,3.825)%
  --(2.050,4.245)--(2.235,4.507)--(2.419,4.532)--(2.604,4.274)--(2.788,3.753)--(2.973,3.074)%
  --(3.157,2.383)--(3.342,1.799)--(3.526,1.377)--(3.711,1.127)--(3.895,1.037)--(4.080,1.097)%
  --(4.264,1.304)--(4.449,1.659)--(4.633,2.162)--(4.818,2.806)--(5.002,3.551)--(5.187,4.309)%
  --(5.371,4.926)--(5.556,5.226)--(5.740,5.101)--(5.925,4.549)--(6.109,3.695)--(6.294,2.766)%
  --(6.478,1.973)--(6.663,1.413)--(6.847,1.096)--(7.032,0.989)--(7.216,1.054)--(7.401,1.259)%
  --(7.585,1.584);
\gppoint{gp mark 9}{(1.497,2.826)}
\gppoint{gp mark 9}{(2.235,4.507)}
\gppoint{gp mark 9}{(2.973,3.074)}
\gppoint{gp mark 9}{(3.711,1.127)}
\gppoint{gp mark 9}{(4.449,1.659)}
\gppoint{gp mark 9}{(5.187,4.309)}
\gppoint{gp mark 9}{(5.925,4.549)}
\gppoint{gp mark 9}{(6.663,1.413)}
\gppoint{gp mark 9}{(7.401,1.259)}
\gppoint{gp mark 9}{(6.943,5.779)}
\gpcolor{color=gp lt color border}
\gpsetlinewidth{1.00}
\draw[gp path] (1.128,6.306)--(1.128,0.739)--(7.585,0.739)--(7.585,6.306)--cycle;
\gpdefrectangularnode{gp plot 1}{\pgfpoint{1.128cm}{0.739cm}}{\pgfpoint{7.585cm}{6.306cm}}
\end{tikzpicture}
  \end{minipage}
  \begin{minipage}[b]{0.47\textwidth}
    \input{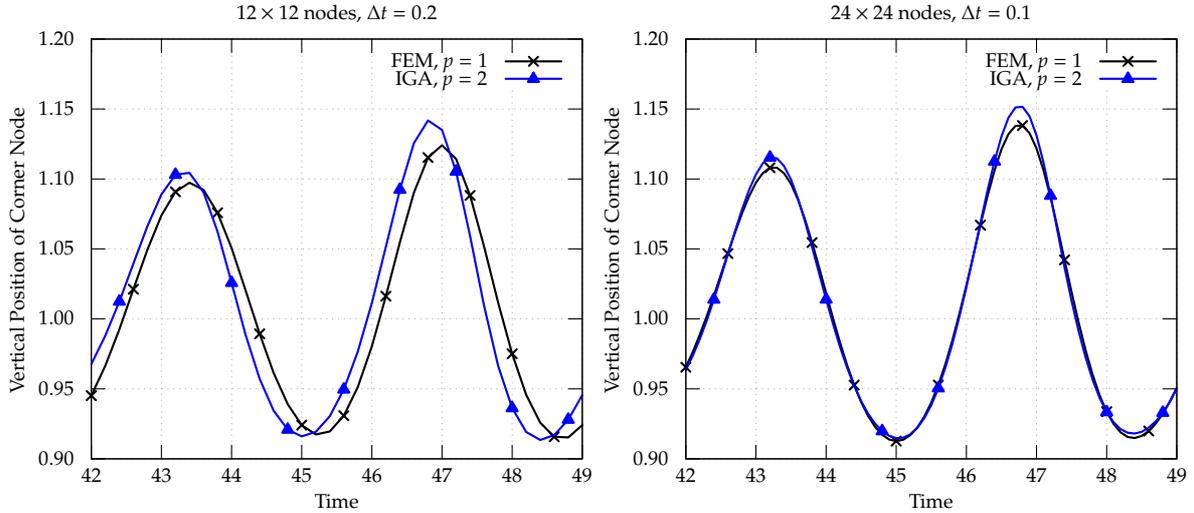}
  \end{minipage}
  \caption{Comparison of linear FEM and quadratic order IGA simulations with
           PDE-based normal-direction movement on
           coarse mesh (left) and
           fine mesh (right).
           The meshes have been refined equally in spatial and temporal
           dimensions in order to keep the element size ratio roughly constant.}
  \label{fig:tank.femvsiga.nodemov}
\end{figure}

\begin{figure}[p]
  \centering
  \begin{minipage}[b]{\textwidth}
    \centering
    \begin{minipage}[b]{0.24\textwidth}
      \centering
      \begin{tikzpicture}[scale=3.3]
        \tikzstyle{mesh} = [thin]
        \input{Meshes/cos.iga.12x12.normal.dt100e-3/mesh1.tex}
      \end{tikzpicture}

      \footnotesize
      $t=1.7$
    \end{minipage}
    \begin{minipage}[b]{0.24\textwidth}
      \centering
      \begin{tikzpicture}[scale=3.3]
        \tikzstyle{mesh} = [thin]
        \input{Meshes/cos.iga.12x12.normal.dt100e-3/mesh2.tex}
      \end{tikzpicture}

      \footnotesize
      $t=18.0$
    \end{minipage}
    \begin{minipage}[b]{0.24\textwidth}
      \centering
      \begin{tikzpicture}[scale=3.3]
        \tikzstyle{mesh} = [thin]
        \input{Meshes/cos.iga.12x12.normal.dt100e-3/mesh3.tex}
      \end{tikzpicture}

      \footnotesize
      $t=23.6$
    \end{minipage}
    \begin{minipage}[b]{0.24\textwidth}
      \centering
      \begin{tikzpicture}[scale=3.3]
        \tikzstyle{mesh} = [thin]
        \input{Meshes/cos.iga.12x12.normal.dt100e-3/mesh4.tex}
      \end{tikzpicture}

      \footnotesize
      $t=49.9$
    \end{minipage}
    \caption{Element quality in physical space with Greville-based normal
direction movement ($12 \times 12$ control points, $\Delta t = 0.1$). The
surface shows a kink.}
    \label{fig:tank.meshes.greville}
  \end{minipage}

  \vspace{.5cm}

  \begin{minipage}[b]{\textwidth}
    \centering
    \begin{minipage}[b]{0.24\textwidth}
      \centering
      \begin{tikzpicture}[scale=3.3]
        \tikzstyle{mesh} = [thin]
        \input{Meshes/cos.iga.12x12.pdenor.dt100e-3/mesh1.tex}
      \end{tikzpicture}

      \footnotesize
      $t=1.7$
    \end{minipage}
    \begin{minipage}[b]{0.24\textwidth}
      \centering
      \begin{tikzpicture}[scale=3.3]
        \tikzstyle{mesh} = [thin]
        \input{Meshes/cos.iga.12x12.pdenor.dt100e-3/mesh2.tex}
      \end{tikzpicture}

      \footnotesize
      $t=18.0$
    \end{minipage}
    \begin{minipage}[b]{0.24\textwidth}
      \centering
      \begin{tikzpicture}[scale=3.3]
        \tikzstyle{mesh} = [thin]
        \input{Meshes/cos.iga.12x12.pdenor.dt100e-3/mesh3.tex}
      \end{tikzpicture}

      \footnotesize
      $t=23.6$
    \end{minipage}
    \begin{minipage}[b]{0.24\textwidth}
      \centering
      \begin{tikzpicture}[scale=3.3]
        \tikzstyle{mesh} = [thin]
        \input{Meshes/cos.iga.12x12.pdenor.dt100e-3/mesh4.tex}
      \end{tikzpicture}

      \footnotesize
      $t=49.9$
    \end{minipage}
    \caption{Element quality in physical space with PDE-based normal direction
movement. The elements become distorted over time.  ($12 \times 12$ control
points, $\Delta t = 0.1$).}
    \label{fig:tank.meshes.normal}
  \end{minipage}

  \vspace{.5cm}

  \begin{minipage}[b]{\textwidth}
    \centering
    \begin{minipage}[b]{0.24\textwidth}
      \centering
      \begin{tikzpicture}[scale=3.3]
        \tikzstyle{mesh} = [thin]
        \input{Meshes/cos.iga.12x12.pdedir.dt100e-3/mesh1.tex}
      \end{tikzpicture}

      \footnotesize
      $t=1.7$
    \end{minipage}
    \begin{minipage}[b]{0.24\textwidth}
      \centering
      \begin{tikzpicture}[scale=3.3]
        \tikzstyle{mesh} = [thin]
        \input{Meshes/cos.iga.12x12.pdedir.dt100e-3/mesh2.tex}
      \end{tikzpicture}

      \footnotesize
      $t=18.0$
    \end{minipage}
    \begin{minipage}[b]{0.24\textwidth}
      \centering
      \begin{tikzpicture}[scale=3.3]
        \tikzstyle{mesh} = [thin]
        \input{Meshes/cos.iga.12x12.pdedir.dt100e-3/mesh3.tex}
      \end{tikzpicture}

      \footnotesize
      $t=23.6$
    \end{minipage}
    \begin{minipage}[b]{0.24\textwidth}
      \centering
      \begin{tikzpicture}[scale=3.3]
        \tikzstyle{mesh} = [thin]
        \input{Meshes/cos.iga.12x12.pdedir.dt100e-3/mesh4.tex}
      \end{tikzpicture}

      \footnotesize
      $t=49.9$
    \end{minipage}
    \caption{Element quality in physical space with PDE-based directional
(vertical) movement.  ($12 \times 12$ control points, $\Delta t = 0.1$). The
element quality remains high throughout the entire simulation.}
    \label{fig:tank.meshes.directional}
  \end{minipage}

  \vspace{.5cm}

  \begin{minipage}[b]{\textwidth}
    \centering
    \begin{minipage}[b]{0.24\textwidth}
      \centering
      \begin{tikzpicture}[scale=3.3]
        \tikzstyle{mesh} = [thin]
        \input{Meshes/cos.iga.12x12.pdeexa.dt100e-3/mesh1.tex}
      \end{tikzpicture}

      \footnotesize
      $t=1.7$
    \end{minipage}
    \begin{minipage}[b]{0.24\textwidth}
      \centering
      \begin{tikzpicture}[scale=3.3]
        \tikzstyle{mesh} = [thin]
        \input{Meshes/cos.iga.12x12.pdeexa.dt100e-3/mesh2.tex}
      \end{tikzpicture}

      \footnotesize
      $t=18.0$
    \end{minipage}
    \begin{minipage}[b]{0.24\textwidth}
      \centering
      \begin{tikzpicture}[scale=3.3]
        \tikzstyle{mesh} = [thin]
        \input{Meshes/cos.iga.12x12.pdeexa.dt100e-3/mesh3.tex}
      \end{tikzpicture}

      \footnotesize
      $t=23.6$
    \end{minipage}
    \begin{minipage}[b]{0.24\textwidth}
      \hspace{1cm}
    \end{minipage}
    \caption{Element quality in physical space with PDE-based equal movement.
($12 \times 12$ control points, $\Delta t = 0.1$).The elements become distorted
over time. }
    \label{fig:tank.meshes.equal}
  \end{minipage}

\end{figure}

In addition to mass conservation, the overall shape of the free surface is of
interest. One possible indicator is the vertical position --- $y(t)$ --- of the
northeast corner of the simulation domain plotted over time (cf.
Figure~\ref{fig:tank.mov} for the results on the fine meshes).
For standard FEM, all versions lead to qualitatively similar results (top-left
plot). Only a zoom (lower-left plot) reveals minor differences.
For IGA --- as for FEM --- all PDE versions lead to similar results (top-right
plot), where again only a zoom reveals minor variations (lower-right plot).
However, the Greville version yields results that are in no relation to all
other results.

So far, standard FEM and IGA were always regarded separately, since there are
plenty of reasons for differing results --- e.g. degree of interpolation,
continuity across elements, imposition of boundary conditions --- that go far
beyond the task of mesh adaptation. Nevertheless --- in the limit of fine
temporal and spatial resolution ---, the results for all variants are expected
to converge towards one common solution.
Since the exact solution is not available for comparison, it is difficult to
judge the correctness of the respective results. Still, it is possible to
observe the results during mesh refinement.
Figure~\ref{fig:tank.femvsiga.nodemov} shows the results for standard FEM and
IGA with PDE-based normal-direction movement on the respective coarse and fine
meshes. On the coarse meshes, the results show significant differences in both
phase and amplitude. While it is impossible to say which variant is closer to
the exact solution, the fact that the results are far more in agreement on the
fine meshes indicates that both variants converge towards the same limit.

When using the total mass flux as a criterion for judging the quality of the
method, both spatially and temporally local errors in the no penetration
condition can be overlooked. For this reason, we introduce another error
measure: the flux error $\tilde{f}$ over the space-time boundary. This is
defined as

\begin{equation}
  \tilde{f} = \sqrt{\frac{\int_0^T \int_{\Gamma_{free,t}}
              \left| {\bf u} \cdot {\bf n} - {\bf v} \cdot {\bf n} \right|^2
              \, \text{d}{\bf x} \, \text{d}t}
              {\int_0^T \int_{\Gamma_{free,t}} \text{d}{\bf x} \, \text{d}t}} \,.
\end{equation}

For the example of an IGA simulation with directional movement, we show
the flux error $\tilde{f}$ for different mesh sizes and time step widths
in Table~\ref{fig:tank.timesteps}. The table shows an almost perfect linear
relationship between the flux error and the time step width. The influence of
the mesh resolution only becomes relevant for extremely small time step widths,
which are no longer practical.
One should note that this error measure only considers the flux over the free
surface, and not the quality of the overall solution. The finer mesh can still
be expected to be beneficial for the flow solution.

\begin{table}[h!]
  \centering
  \begin{tabular}{| c | r | r | r | r | r | r |}
    \hline
    Mesh \textbackslash \, $\Delta t$
    & $2.0  \cdot 10^{-1}$
    & $1.0  \cdot 10^{-1}$
    & $ 5.0  \cdot 10^{-2}$
    & $ 2.5  \cdot 10^{-2}$
    & $ 1.25  \cdot 10^{-2}$
    & $  6.25 \cdot 10^{-3}$ \\
    \hline
    $12 \times 12$
    & $0.009404$
    & $0.004748$
    & $0.002393$
    & $0.001229$
    & $0.000687$
    & $0.000482$ \\
    \hline
    $24 \times 24$
    & $0.009399$
    & $0.004743$
    & $0.002377$
    & $0.001191$
    & $0.000599$
    & $0.000307$ \\
    \hline
  \end{tabular}
  \caption{Table showing flux error $\tilde{f}$ for different time step widths and mesh
           resolutions for an IGA simulation of the sloshing tank with
           directional free-surface movement.}
  \label{fig:tank.timesteps}
\end{table}

With mass conservation and final result being in very similar ranges, up until
now, it remains unclear which PDE-based mesh deformation version is the most
suited for the application of sloshing tank. An answer can be found when
comparing the mesh quality for the different versions. Here, the IGA variants
will serve as an example, but all effects can be observed in the standard FEM
meshes in the same way.
Figures~\ref{fig:tank.meshes.greville} through \ref{fig:tank.meshes.equal} show
clear advantages of the PDE-based vertical displacement over all other versions.
As already seen before, the Greville-based displacement leads to a free surface
with a small kink, which is considered to be an invalid result
(Figure~\ref{fig:tank.meshes.greville}). PDE-based normal
(Figure~\ref{fig:tank.meshes.normal}) and PDE-based equal
(Figure~\ref{fig:tank.meshes.equal}) displacement both have the disadvantage
that they --- to some extent --- still allow for tangential displacement of
nodes / control points. Over the course of several sloshing cycles, this
inevitably leads to poor mesh quality --- note the distorted elements --- or
even breakdown of the simulation. Only the PDE-based vertical displacement leads
to a mesh displacement that can be carried out with no practical limit. All
results were plotted for the coarse mesh for the sake of better readability, but
hold in the same way for the fine mesh.

\subsection{Test case die swell}

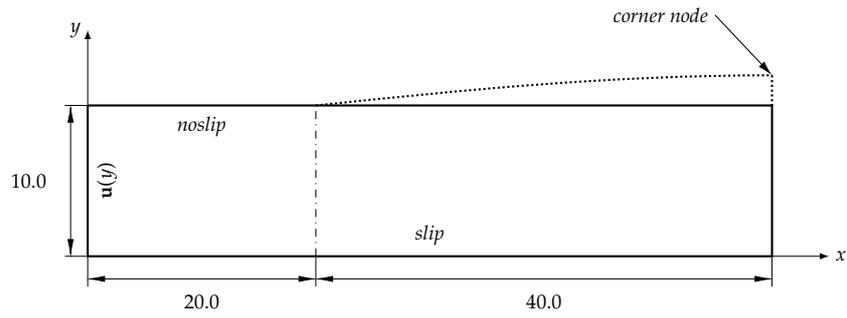
\begin{figure}[h!]
  \centering
  \begin{tikzpicture}[xscale=0.15, yscale=0.2]

    \tikzset{every node/.append style={scale=0.75}}

    \draw [thick] (0.0,10.0) -- (0.0,0.0) -- (60.0,0.0) -- (60.0,10.0) -- (0.0,10.0);
    \draw [thin, dash dot] (20.0,0.0) -- (20.0,10.0);

    \draw [thick, densely dotted] (20,10) sin (60,12);
    \draw [thick, densely dotted] (60,12) -- (60,10);

    \draw [very thin] (0.0,0.0) -- (0.0, -2.0);
    \draw [very thin] (20.0,0.0) -- (20.0, -2.0);
    \draw [{Latex[length=2.25mm,width=0.75mm]}-{Latex[length=2.25mm, width=0.75mm]}, very thin] (0.0,-1.5) -- (20.0,-1.5);
    \node [thick] at (10.0, -3.0) {$20.0$};
    \draw [very thin] (60.0,0.0) -- (60.0, -2.0);
    \draw [{Latex[length=2.25mm,width=0.75mm]}-{Latex[length=2.25mm, width=0.75mm]}, very thin] (20.0,-1.5) -- (60.0,-1.5);
    \node [thick] at (40.0, -3.0) {$40.0$};
    \draw [very thin] (0.0,0.0) -- (-2.0, 0.0);
    \draw [very thin] (0.0,10.0) -- (-2.0, 10.0);
    \draw [{Latex[length=2.25mm,width=0.75mm]}-{Latex[length=2.25mm, width=0.75mm]}, very thin] (-1.5,0.0) -- (-1.5,10);
    \node [thick, left] at (-3.0,5.0) {$10.0$};

    \draw [very thin, >=latex, ->] (60.0,0.0) -- (65.0,0.0) node[right]{$x$};
    \draw [very thin, >=latex, ->] (0.0,10.0) -- (0.0,15.0) node[left]{$y$};
    \node [thick, rotate=90] at (1.7,5.0) {${\bf u}(y)$};
    \node [thick] at (10.0,8.6) {$noslip$};
    \node [thick] at (30.0,1.4) {$slip$};
    \draw [{Latex[length=2mm, width=1mm]}-, thin] (60,12) -- (55,16) node[left]{$corner\ node$};

  \end{tikzpicture}
  \caption{Die swell: test case setup. Velocity on left boundary
  ${\bf u} \left( y \right) = \left(
0.1\left(100-y^2\right),0\right)^T$. Schematic quasi-steady solution shown by
dotted line.}
  \label{fig:swell.testcase}
\end{figure}

The second test case considers the behavior of a fluid that exits a nozzle (cf.
Figure~\ref{fig:swell.testcase} for the exact set-up). The Newtonian fluid ---
density 1.0 and viscosity 100000 --- enters the nozzle from the left with a
prescribed velocity of ${\bf u}(y)=\left(0.1\left(100-y^2\right),0\right)^T$.
Gravity is neglected. The lower boundary is a symmetry boundary with a slip
condition on the velocity. The top boundary is partially a wall --- nozzle part
--- and partially a free surface.
On the free-surface part, natural boundary conditions are used. This means that
the stress on the boundary will be zero in surface normal direction.
Due to the change in boundary condition on the
top boundary from no-slip to free surface, the velocity profile restructures
itself from parabolic --- pipe flow --- to a block profile. This introduces
velocity components orthogonal to the flow direction, which in turn leads to the
so-called die swell. Note that the orthogonal velocity component is orders of
magnitude smaller than the velocity in the main flow direction.

\begin{figure}[h!]
  \centering
  \begin{minipage}[b]{\textwidth}
    \centering
    \includegraphics{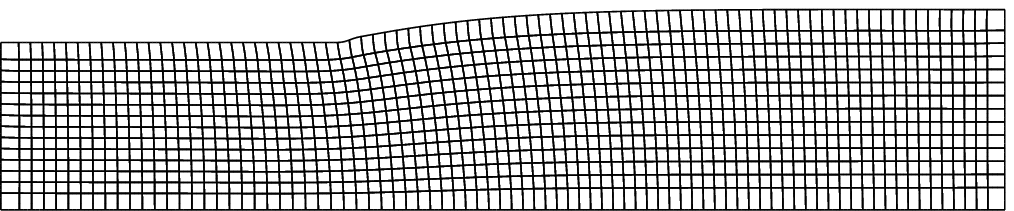}

    \footnotesize
    $86 \times 16$ PDE-based directional (vertical) movement at $t=14.25$
  \end{minipage}

  \vspace{.5cm}

  \begin{minipage}[b]{\textwidth}
    \centering
    \includegraphics{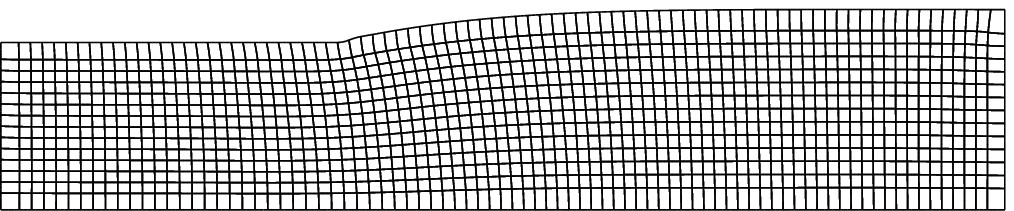}

    \footnotesize
    $86 \times 16$ PDE-based normal movement at $t=14.25$
  \end{minipage}
  \caption{Meshes used for IGA simulations with $p=2$. Elements are shown in physical
           space. While both meshes retain good quality, the directional movement seems
           to yield a better element shape in the upper right corner.}
  \label{fig:swell.meshes}
\end{figure}

Again, standard FEM and IGA computations were performed. For standard FEM, we
utilized a Q1Q1 structured mesh with $86$ by $16$ nodes. For IGA, a NURBS of
degree $2$ with $86$ by $16$ control points was employed. In both cases, the
time step was $0.015625$. Note that in both cases, PDE-based displacement with
the full velocity is not shown, as it is not usable for this application due to
the strong tangential velocity component.
The final meshes are given in Figure~\ref{fig:swell.meshes}.
Clearly, either of the proposed methods suitable in this case, retains good mesh
quality.

\begin{figure}[h!]
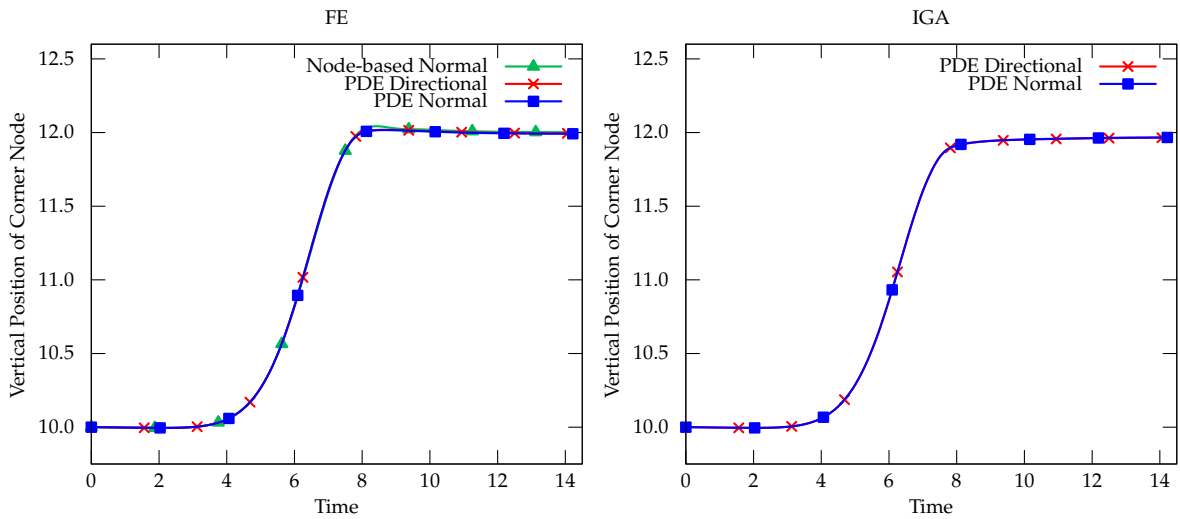

  \centering
  \begin{minipage}[b]{0.47\textwidth}
    \input{Plots/Swell/fem.86x16.dt0.015625.nodeMov.log.tex}
  \end{minipage}
  \begin{minipage}[b]{0.47\textwidth}
    \input{Plots/Swell/iga.86x16.dt0.015625.nodeMov.log.tex}
  \end{minipage}
  \caption{Comparison of node movement for FE and IGA}
  \label{fig:swell:femvsiga.mov}
\end{figure}

Figure~\ref{fig:swell:femvsiga.mov} plots the position of the corner node --- as
indicated in Figure~\ref{fig:swell.testcase} --- over time. On the left,
standard FEM results are displayed, on the right, the IGA results.
All applied methods lead to nearly identical results.
Thus, for this test case, no clear recommendation for any particular method can
be given. However --- as for the sloshing tank --- the introduction of the
PDE-based displacement methods allow for the usage of IGA.


\section{Conclusion: Which mesh displacement method should be used?} \label{sec-conc}

The paper addresses options for interface tracking of free-surface problems. When considering boundary-conforming meshes in such a context, the adaptation of the boundary position --- i.e., the boundary deformation --- based on the flow field, or --- more specifically --- the no-penetration boundary condition, is the major challenge. Motivated by the integral formulation of the no-penetration boundary condition, a PDE-based boundary displacement method is proposed. The new method has two significant advantages over the current state of the art: (1) It is applicable to finite element methods of linear and higher order as well as to isogeometric finite element methods, and (2) It guarantees mass conservation up to numerical errors caused by the solution process of the non-linear equation system, i.e., usually up to machine precision. In two space dimensions, the validity of the presented method has been confirmed with two standard numerical test cases for free-surface flow. 
The basic idea of the method is also applicable to three-dimensional problems.

\section*{Acknowledgments}

The authors gratefully acknowledge the support of DFG under the Collaborative
Research Center SFB 1120 (subproject B2) and DFG grant "Automated design and
optimisation of dynamic mixing and shear elements for single-screw extruders"
(EL 741/5-1).

\bibliographystyle{elsarticle-num}
\bibliography{references}

\begin{thebibliography}{10}
\expandafter\ifx\csname url\endcsname\relax
  \def\url#1{\texttt{#1}}\fi
\expandafter\ifx\csname urlprefix\endcsname\relax\def\urlprefix{URL }\fi
\expandafter\ifx\csname href\endcsname\relax
  \def\href#1#2{#2} \def\path#1{#1}\fi

\bibitem{Elgeti2015}
S.~Elgeti, H.~Sauerland, {Deforming Fluid Domains Within the Finite Element
  Method: Five Mesh-Based Tracking Methods in Comparison}, Archives of
  Computational Methods in Engineering 23~(2) (2015) 323--361.

\bibitem{Caboussat2005}
A.~Caboussat, Numerical simulation of two-phase free surface flows 12~(2)
  (2005) 165--24.

\bibitem{Easton72}
C.~R. Easton, {Homogeneous Boundary Conditions for Pressure in the MAC Method}
  9 (1972) 375--379.

\bibitem{Girault76}
V.~Girault, {A Combined Finite Element and Marker and Cell Method for Solving
  Navier-Stokes Equations} 26 (1976) 39--59.

\bibitem{Sethian99b}
J.~Sethian, Level Set Methods and Fast Marching Methods, 2nd Edition, Cambridge
  University Press, 1999.

\bibitem{Osher2001}
S.~Osher, R.~Fedkiw, Level set methods: An overview and some recent results 169
  (2001) 463--502.

\bibitem{Nichols71}
B.~Nichols, C.~Hirt, Improved free surface boundary conditions for numerical
  incompressible-flow calculations 8~(3) (1971) 434--448.

\bibitem{Hirt75}
C.~W. Hirt, B.~D. Nichols, N.~C. Romero, {SOLA}: A numerical solution algorithm
  for transient fluid flows, Tech. Rep. 32418, {NASA STI/}Recon Technical
  Report N 75 (1975).

\bibitem{Noh76}
W.~Noh, P.~Woodward, {SLIC (Simple Line Interface Calculation)}, Lecture Notes
  in Physics 59 (1976) 330--340.

\bibitem{Hirt81a}
C.~Hirt, B.~Nichols, Volume of fluid ({VOF}) method for the dynamics of free
  boundaries 39 (1981) 201 -- 225.

\bibitem{Gibbs1873}
J.~Gibbs, {On the Equilibrium of Heterogeneous Substances}, {Transactions of
  the Conneticut Academy of Arts and Sciences}, 1874-1878.

\bibitem{Onsager1931}
L.~Onsager, {Reciprocal Relations in Irreversible Processes}, Physical Review
  37~(405).

\bibitem{Prigogine1966}
I.~Prigogine, {Non-Equilibrium Statistical Mechanics}, 2nd Edition, Wiley, New
  York, 1966.

\bibitem{Gyrmati1970}
I.~Gyrmati, {Non-Equilibrium Thermodynamics}, 1st Edition, Springer, New York,
  1970.

\bibitem{deGroot1984}
S.~de~Groot, P.~Mazur, {Non-Equilibrium Thermodynamics}, 1st Edition, Dover,
  New York, 1984.

\bibitem{Emmerich2002}
H.~Emmerich, {The Diffuse Interface Approach in Matrials Science: Thermodynamic
  Concepts and Applications of Phase-Field Models}, Springer, 2003.

\bibitem{Hirt74a}
C.~Hirt, A.~Amsden, J.~Cook, An arbitrary {L}agrangian {E}ulerian computing
  method for all flow speeds 14 (1974) 227 -- 253.

\bibitem{Tezduyar92a}
T.~E. Tezduyar, M.~Behr, J.~Liou, A new strategy for finite element
  computations involving moving boundaries and interfaces--the
  deforming-spatial-domain/space-time procedure: {I}. the concept and the
  preliminary numerical tests 94~(3) (1992) 339 -- 351.

\bibitem{Grotle2016}
E.~Grotle, H.~Bihs, E.~P. an~V.~Aesoy, {CFD Simulations of Non-Linear Sloshing
  in a Rotating Rectangular Tank Using the Level Set Method}, in: ASME 2016
  35th International Conference on Ocean, Offshore and Arctic Engineering,
  2016.

\bibitem{Manderbacka2016}
J.~Fonfach, T.~Manderbacka, M.~Neves, {Numerical sloshing simulations:
  Comparison between lagrangian and lumped mass models applied to two
  compartments with mass transfer}, Ocean Engineering 114 (2016) 168--184.

\bibitem{Elgeti2010}
S.~Elgeti, H.~Sauerland, L.~Pauli, M.~Behr, On the usage of {NURBS} as
  interface representation in free-surface flows 69~(1) (2012) 73--87.

\bibitem{Pauli2012}
L.~Pauli, M.~Behr, S.~Elgeti, Towards shape optimization of extrusion dies with
  respect to homogeneous die swell, {Journal of Non-Newtonian Fluid Mechanics}
  69 (2012) 73--87.

\bibitem{Stavrev2015}
A.~Stavrev, P.~Knechtges, S.~Elgeti, A.~Huerta, Space-time nurbs-enhanced
  finite elements for free-surface flows in 2d, {International Journal for
  Numerical Methods in Fluids} 81~(7) (2016) 397--459.

\bibitem{Knechtges2014}
P.~Knechtges, M.~Behr, S.~Elgeti, {Fully-implicit log-conformation formulation
  of constitutive laws}, {Journal of Non-Newtonian Fluid Mechanics} 214 (2014)
  78--87.

\bibitem{Knechtges2015}
P.~Knechtges, {The fully-implicit log-conformation formulation and its
  application to three-dimensional flows}, {Journal of Non-Newtonian Fluid
  Mechanics} 223 (2015) 209--220.

\bibitem{Tome2008}
M.~Tom{\'e}, A.~Castelo, V.~Ferreira, S.~McKee, {A finite difference technique
  for solving the Oldroyd-B model for 3D-unsteady free surface flows}, Journal
  of Non-Newtonian Fluid Mechanics 154 (2008) 179--206.

\bibitem{Tome2012}
M.~Tom{\'e}, A.~Castelo, A.~Afonso, M.~Alves, F.~Pinho, Application of the
  log-conformation tensor to three-dimensional time-dependent free surface
  flows, Journal of Non-Newtonian Fluid Mechanics.

\bibitem{Mompean2011}
G.~Mompean, L.~Thais, M.~Tom{\'e}, A.~Castelo, Numerical prediction of
  three-dimensional time-dependent viscoelastic extrudate swell using
  differential and algebraic models, Computers \& Fluids 44 (2011) 68--78.

\bibitem{Rocco2010}
G.~Rocco, G.~Coppola, L.~de~Luca, {The VOF method applied to the numerical
  simulation of a 2D liquid jet under gravity}, {WIT Transactions on
  Engineering Sciences} 69.

\bibitem{Hughes2005}
T.~J.~R. Hughes, J.~A. Cottrell, Y.~Bazilevs, {Isogeometric analysis: CAD,
  finite elements, NURBS, exact geometry and mesh refinement} 194 (2005)
  4135--4195.

\bibitem{Cottrell2009}
J.~A. Cottrell, T.~J.~R. Hughes, Y.~Bazilevs, {Isogeometric Analysis: Toward
  Integration of CAD and FEA}, {John Wiley \& Sons, Ltd}, 2009.

\bibitem{Gomez2008}
H.~G{\'o}mez, V.~Calo, Y.~Bazilevs, T.~Hughes, {Isogeometric analysis of the
  Cahn--Hilliard phase-field model} 197~(49) (2008) 4333--4352.

\bibitem{Borden2012}
M.~Borden, J.~Michael, C.~Verhoosel, M.~Scott, T.~J. Hughes, C.~Landis, {A
  phase-field description of dynamic brittle fracture}, {Computer Methods in
  Applied Mechanics and Engineering} 217 (2012) 77--95.

\bibitem{Dede2012}
L.~Ded{\`e}, M.~Borden, T.~Hughes, {Isogeometric Analysis for Topology
  Optimization with a Phase Field Model} 19 (2012) 427--465.

\bibitem{Akkerman2011}
I.~Akkerman, Y.~Bazilevs, C.~E. Kees, M.~W. Farthing, Isogeometric analysis of
  free-surface flow, {Journal of Computational Physics} 11 (2011) 4137--4152.

\bibitem{Amini2014}
R.~Amini, R.~Maghsoodi, N.~Z. Moghaddamog, {Simulating free surface problem
  using isogeometric analysis}, {Journal of the Brazilian Society of Mechanical
  Sciences and Engineering} 38~(2) (2016) 413--421.

\bibitem{Rank2012}
E.~Rank, M.~Ruess, S.~Kollmannsberger, D.~Schillinger, A.~D\"uster, Geometric
  modeling, isogeometric analysis and the finite cell method, {Computer Methods
  in Applied Mechanics and Engineering} 249 (2012) 104--115.

\bibitem{Ginnis2014}
A.~Ginnisa, K.~Kostasb, C.~Politisb, P.~Kaklisa, K.~Belibassakisa,
  T.~Gerostathisb, M.~Scott, T.~Hughes, {Isogeometric Boundary-Element Analysis
  for the Wave-Resistance Problem using T-splines}, {Computer Methods in
  Applied Mechanics and Engineering} 279 (2014) 425--439.

\bibitem{Behr92}
M.~Behr, {Stabilized Finite Element Methods for Incompressible Flows with
  Emphasis on Moving Boundaries and Interfaces}, Ph.D. thesis, University of
  Minnesota, Department of Aerospace Engineering and Mechanics (1992).

\bibitem{Johnson94a}
A.~Johnson, T.~Tezduyar, Mesh update strategies in parallel finite element
  computations of flow problems with moving boundaries and interfaces 119
  (1994) 73 -- 94.

\bibitem{Tezduyar92b}
T.~E. Tezduyar, M.~Behr, S.~Mittal, J.~Liou, A new strategy for finite element
  computations involving moving boundaries and interfaces--the
  deforming-spatial-domain/space-time procedure: {II}. computation of
  free-surface flows, two-liquid flows, and flows with drifting cylinders
  94~(3) (1992) 353 -- 371.

\bibitem{Tezduyar2004}
T.~Tezduyar, Computation of moving boundaries and interfaces and stabilization
  parameters, {International Journal for Numerical Methods in Fluids} 43 (2004)
  555--575.

\bibitem{Takizawa2011}
K.~Takizawa, B.~Henicke, T.~Tezduyar, M.-C. Hsu, Y.~Bazilevs, {Stabilized
  Space--Time Computation of Wind-Turbine Rotor Aerodynamics}, Computational
  Mechanics 48 (2011) 333--344.

\bibitem{Takizawa2014}
K.~Takizawa, T.~Tezduyar, A.~Buscher, S.~Asada, {Space-time interface tracking
  with topology change}, Computational Mechanics 54 (2014) 955--971.

\bibitem{Zilian2008}
A.~Zilian, A.~Legay, The enriched space--time finite element method (est) for
  simultaneous solution of fluid--structure interaction, International Journal
  for Numerical Methods in Biomedical Engineering 75~(3).

\bibitem{Engelman82a}
M.~Engelman, R.~Sani, P.~Gresho, The implementation of normal and/or tangential
  boundary conditions in finite element codes for incompressible fluid flow 2
  (1982) 225--238.

\bibitem{Piegl97}
L.~Piegel, W.~Tiller, The NURBS Book, Springer, Berlin, Germany, 1997.

\bibitem{Rogers2001}
D.~Rogers, An Introduction to NURBS with Historical Perspective, Morgan
  Kaufmann Publishers, 2001.

\bibitem{Bazilevs2010}
Y.~Bazilevs, V.~Calo, J.~Cottrell, J.~Evans, T.~Hughes, S.~Lipton, M.~Scott,
  T.~Sederberg, {Isogeometric analysis using T-splines} 199 (2010) 229--263.

\bibitem{Bazilevs2007}
Y.~Bazilevs, T.~J. Hughes, {Weak imposition of Dirichlet boundary conditions in
  fluid mechanics}, {Computers and Fluids} 36 (2007) 12--26.

\bibitem{Farin2014}
G.~Farin, {Curves and Surfaces for Computer-Aided Geometric Design: A Practical
  Guide}, Elsevier, 2014.

\bibitem{Gray1998}
A.~Gray, Modern differential geometry of curves and surfaces with Mathematica,
  2nd Edition, CRC Press, Boca Raton, 1998.

\end{thebibliography}


\end{document}